\documentclass[mnsc]{informsarxiv}
\DeclareFontShape{OT1}{cmr}{m}{n}{<-6>cmr5<6-7>cmr6<7-8>cmr7<8-9>cmr8<9-10>cmr9<10->cmr10}{}
\usepackage[authoryear,round]{natbib}
\usepackage[margin=1in]{geometry}
\usepackage{mathtools}
\usepackage{algorithm,algpseudocode}
\usepackage{caption,subcaption}
\usepackage{enumitem}
\usepackage{prettyref}
\usepackage[bookmarks,colorlinks=true,plainpages=false,hypertexnames=false,citecolor=blue,linkcolor=blue,urlcolor=blue]{hyperref}
\theoremstyle{plain}
\newtheorem{theorem}{Theorem}
\newtheorem{lemma}{Lemma}
\newtheorem{proposition}{Proposition}
\newtheorem{corollary}{Corollary}
\theoremstyle{definition}
\newtheorem{definition}{Definition}

\newcommand{\E}{\mathbb{E}}
\newcommand{\reals}{{\mathbb{R}}}
\newcommand{\diff}{{\rm d}}
\newcommand{\expect}[1]{\mathbb{E}\left[ #1 \right]}
\newcommand{\norm}[1]{\left\|{#1} \right\|}
\newcommand{\cnorm}[1]{\left\| #1 \right\|_{\mathcal{I}^+}}
\newcommand{\iprod}[2]{\left \langle #1, #2 \right\rangle}
\newcommand{\calB}{{\mathcal{B}}}
\newcommand{\calI}{{\mathcal{I}}}
\newcommand{\dualopt}{U^*}
\renewcommand{\tilde}{\widetilde}
\renewcommand{\bar}{\overline}
\renewcommand{\hat}{\widehat}

\newrefformat{eq}{(\ref{#1})}
\newrefformat{sec}{Section~\ref{#1}}
\newrefformat{lmm}{Lemma~\ref{#1}}

\RUNTITLE{Dynamic Matching with Abandonments}
\TITLE{Longest-Queue Policies in Dynamic Matching with Abandonments} 

\ARTICLEAUTHORS{
\AUTHOR{Siqi He, Yehua Wei, Jiaming Xu}
\AFF{The Fuqua School of Business, Duke University, 100 Fuqua Drive, Durham, North Carolina 27708,\\
\EMAIL{siqi.he@duke.edu}, \EMAIL{yehua.wei@duke.edu}, \EMAIL{jiaming.xu868@duke.edu}}	
}

\ABSTRACT{We study a dynamic two-way matching market in which agents of heterogeneous types
arrive over time, wait to be matched according to a compatibility network,
and abandon the system if their random patience expires. 
The decision maker
collects match rewards and pays holding and abandonment costs, and seeks to
maximize the expected net value over a finite horizon. We evaluate policies
by their approximation ratio against a fluid benchmark, the optimal value
of a linear-programming relaxation that upper bounds the expected net value
of every policy. 
Under a general position gap (GPG) condition, which
requires the unique optimal dual solution of this linear program to remain
optimal under small perturbations of the arrival rates, we propose a
longest-queue policy.
The policy first uses this dual solution to identify a set of good match
types and a set of critical agent types. Upon each arrival, it considers
only good matches and matches the arrival to the critical partner
with the longest queue subject to availabilities.
In the high-volume regime, in which arrival rates are
scaled by \(s\), we prove that the policy attains an approximation ratio of
at least \(1-O(1/s)\) at every scale \(s\geq1\) and every horizon
\(T\geq1\), so the ratio tends to $1$ as $s$ goes to infinity. 
}

\begin{document}
\date{}
\maketitle
\tableofcontents

\section{Introduction}
\label{sec:introduction}

In dynamic matching markets, waiting can create opportunities for more valuable matches, but it can also lead to costly delays and abandonment. This tradeoff arises in a variety of operational settings, including kidney exchange, ride-sharing platforms, and online dating platforms \citep{ozkan2020dynamic, ashlagi2023matching, rios2023improving, kohlenberg2025cost}. 
To understand how to manage this tradeoff, researchers have studied dynamic matching models in which agents of different types arrive over time and a matching network specifies which types of agents can be paired. 
In such models, a decision maker observes the agents currently in the system and decides which matches to make and which agents to leave waiting. Each match generates a reward, and matched agents leave the market. Unmatched agents incur holding costs per unit of time and abandon the market with a penalty when their waiting times exceed their patience times. The objective of the decision maker is to maximize expected total matching rewards net of holding costs and abandonment penalties.

One stream of research on dynamic matching with abandonment focuses on policies with constant-factor guarantees for long-run average performance \citep{collina2020dynamic,aouad2020dynamic,amanihamedani2024improved}.
These guarantees bound the loss relative to optimal performance, but leave a fixed approximation gap. 
Recently, \citet{amanihamedani2025adaptive} made significant progress 
by developing policies that can make this gap arbitrarily small for bipartite matching models with a fixed number of waiting agent types.
However, their model assumes waiting on only one side of the market, while agents on the other side must be matched immediately upon arrival or leave unmatched.
In addition, the policies are also complex, making them difficult to implement in practice. This motivates the search for simple policies with near-optimal performance in markets where agents of all types can wait and abandon.

The high-volume regime offers a promising setting for developing such policies. 
Also known as the large-market regime, 
it has proved useful for analyzing simple policies \citep{akbarpour2020thickness, blanchet2022asymptotically} and for guiding policy design  \citep{ozkan2020dynamic, ata2025dynamic} in various dynamic matching markets.
In the high-volume regime, the arrival rates of all agent types grow proportionally with a market-size parameter, denoted by $s$.
\citet{aveklouris2021matching} study bipartite matching markets in which agents on both sides can wait and abandon. 
They show that simple matching policies can achieve an approximation ratio converging to 1 as the market grows, under a long-run average objective.\footnote{Their asymptotic
guarantee first takes the market-size parameter $s\to\infty$ and then
the time horizon $T\to\infty$.}
Their policies make matches at scheduled review times with review intervals proportional to $s^{-2/3}$, 
allowing the expected number of arrivals between reviews to grow as $s^{1/3}$ as the market grows. Their analysis also accommodates non-exponential patience-time distributions.

In contrast to the long-run asymptotic optimality guarantees for bipartite matching in \citet{aveklouris2021matching}, we establish explicit performance bounds at a given market size $s$ and over a finite operating horizon $T$. 
We consider general two-way matching networks, including those that are non-bipartite, with exponentially distributed patience times.
Under a general position gap (GPG) condition, which requires the fluid LP's unique optimal dual solution to remain optimal under small perturbations of arrival rates, we show that a simple longest-queue policy achieves an approximation ratio of at least $1-C/s$ relative to the optimal policy, for every $s\geq1$ and $T\geq1$, where $C$ is independent of $s$ and $T$.
Without GPG, the same policy remains asymptotically optimal, with a weaker bound of $O(1/\sqrt{s})$ on the approximation gap. 
In addition to the theoretical guarantees, the key difference between our longest-queue policy and the discrete-review policies of \citet{aveklouris2021matching} is that the  longest-queue policy makes matching decisions immediately upon each arrival, avoiding the additional waiting imposed by scheduled reviews and requiring no choice of a review interval. Our numerical experiments demonstrate the practical value of this approach, as the longest-queue policy achieves higher approximation ratios than both discrete-review benchmarks across the tested instances, market sizes, and review intervals.

The longest-queue policy we propose is guided by a fluid linear program that provides an upper bound on achievable performance. The fluid formulation is adapted from the fluid models for overloaded queueing systems with abandonment \citep{whitt2006fluid,atar2010cmu}, as it balances arrival rates against matching and abandonment rates, with exponential patience making abandonment flows linear in queue lengths.
Before matching begins, we use its optimal dual solution to identify a set of \emph{good} match types to use and a set of \emph{critical} agent types to prioritize. These sets remain fixed throughout the operating horizon, so the policy does not need to solve the linear program after the horizon starts. Upon each arrival, it considers available partners
that form a good match. If any critical partners are available, it selects one from the type with the longest queue. Otherwise, it matches with an available noncritical partner that forms a good match, if one
exists. If no good match is available, the arriving agent waits. 
We note that restricting immediate matching to a selected set of compatible pairs is also central to the restricted-greedy policies of \citet{kohlenberg2024quality}, studied in bipartite markets with demand-dependent rewards.

Our analysis builds upon two simple yet key observations. 
First, under exponential patience times, the expected cumulative abandonment is the time integral of expected queue length. Controlling queues therefore controls both holding costs and abandonment penalties. 
Second, after rewriting the objective in terms of adjusted match rewards, complementary slackness and the policy’s restriction to good matches yield an exact regret representation involving only critical types. 
In particular, noncritical queues may be large, and their aggregate abandonment rates may scale with $s$, but their zero dual prices ensure that they do not contribute to the regret. 
The analysis therefore reduces to bounding the expected critical queue length uniformly in time and the market-size parameter $s$. 
Such a bound keeps the expected loss per unit time bounded, while the fluid benchmark grows proportionally with $s$, yielding the $O(1/s)$ relative performance gap.

To establish the desired uniform bound on critical queue lengths, we adapt a standard Lyapunov-drift argument used in recent work establishing constant regret in discrete-time dynamic matching models without abandonment. 
Although abandonment contributes a
nonpositive term to the drift of the quadratic Lyapunov function,
this contribution does not scale with $s$ at a fixed queue size
and hence cannot by itself yield a queue bound uniform in $s$.
Instead, we apply the analysis of \citet{kerimov2023optimality}
to show that, under the GPG condition, the longest-queue policy
generates negative drift of order $s$ whenever the critical-queue
norm exceeds a threshold independent of $s$.
Thus, the longest-queue policy provides the drift needed to control critical queues as the market grows, and abandonment can only strengthen this drift.

The remaining challenge is to convert this negative-drift condition into a bound that is uniform over all times. 
The standard stochastic comparison arguments used in the matching literature, 
e.g., \citet{gupta2024greedy,kerimov2023optimality}, do not directly extend because abandonment intensities depend on the current queues, so systems initialized in different states generally have different departure dynamics. 
Building on the drift method of \citet{hajek1982hitting}, we instead analyze an exponential Lyapunov function for the continuous-time jump process. Although the aggregate abandonment rate may be unbounded, these jumps only decrease the critical queues. 
For these downward jumps, we absorb the quadratic remainder of the exponential increment into its negative first-order contribution. Upward jumps are controlled by the state-independent aggregate arrival rate. This yields an affine drift inequality for the exponential Lyapunov function, and Dynkin’s formula together with Grönwall’s inequality gives the desired critical-queue bound uniformly in both time and $s$.

\section{Model Setup}
\label{sec:model-setup}
We consider a finite set of agent types
\(\mathcal{N}=[n]=\{1,\ldots,n\}\) and a finite set of match types
\(\mathcal{M}=[d]=\{1,\ldots,d\}\).  
We assume that each match type is two-way, that is, \(m\in\mathcal M\)
is determined by two distinct types \(i,j\in[n]\) such that
\(\mathcal N(m) \triangleq \{i,j\}\). Throughout this paper, we use \(m\) and
\((i,j)\) interchangeably.
The match types can be alternatively represented as a network adjacency matrix
\(M\in\{0,1\}^{n\times d}\), where
\(M_{im}=\mathbf{1}\{i\in\mathcal{N}(m)\}\).

The system maintains a separate, initially empty queue for each agent type \(i\in\mathcal{N}\). Type-\(i\) agents arrive and join the corresponding queue
according to a Poisson process with rate \(\lambda_i> 0\), where
\(\lambda=(\lambda_1,\ldots,\lambda_n)^\top\). 
Upon arrival, each type-\(i\) agent draws a
patience time from the exponential distribution with rate \(\theta_i>0\),
independently of all other agents and of the arrival processes, where
\(\theta=(\theta_1,\ldots,\theta_n)^\top\). 
Executing a match of type \(m\) pairs one agent of each type in
\(\mathcal{N}(m)\) and removes both agents from the queues. 
The value of performing a type-\(m\) match is \(r_m\geq 0\). 
If an agent waits longer than its patience time, it abandons the system. A type-\(i\)
agent incurs holding cost \(c_i\) per unit of time while waiting and an
abandonment penalty \(\gamma_i\) if it leaves without being matched.

We consider the class \(\Pi\) of admissible policies, where each
 \(\pi \in \Pi\) is a sequence of (potentially randomized) actions whose action
 at each time \(t\) is measurable with respect to the \(\sigma\)-field generated
 by the history of arrivals, abandonments, and
matching decisions up to time $t$. However, agents' patience times are not observed before abandonments. 

For each \(i\in\mathcal{N}\) and \(t\geq 0\), let \(A_i(t)\) denote the cumulative
number of type-\(i\) arrivals by time \(t\).  A match type
\(m\in\mathcal{M}\) is \emph{available} if the queues of both agent types in
\(\mathcal{N}(m)\) are nonempty.  
At each agent arrival time \(t\), policy
\(\pi\) uses the observed history to select a (possibly empty) collection of
available matches.  For each \(m\in\mathcal{M}\), let \(X_m^\pi(t)\) denote the
number of times match type \(m\) has been performed by time \(t\) under policy
\(\pi\).  
Let \(B_i^\pi(t)\) denote the cumulative number of
type-\(i\) abandonments by time \(t\) under policy \(\pi\).

Let
\(q_i^\pi(t)\) denote the number of waiting type-\(i\) agents
at time \(t\) under policy \(\pi\).  All counting processes start at zero, and
\(q_i^\pi(0)=0\) for every
\(i\in\mathcal{N}\).  The queue length of each type satisfies
\begin{equation}
  q_i^\pi(t)
  =A_i(t)-\sum_{m=1}^d M_{im}X_m^\pi(t)-B_i^\pi(t),
  \qquad i\in\mathcal{N},\quad t\geq 0.
  \label{eq:flow-balance}
\end{equation}

For a queue vector \(q\in\mathbb Z_+^n\), define the \emph{terminal
inventory cost}
\begin{equation*}
    C^{\mathrm{term}}(q)
    \coloneqq
    \sum_{i=1}^n
    \left(\frac{c_i}{\theta_i}+\gamma_i\right)q_i.
\end{equation*}
This cost evaluates the system under the convention that matching stops at
the terminal time and every remaining agent subsequently waits until
abandonment.  By the memoryless property of the exponential patience-time
distribution, a remaining type-\(i\) agent has expected residual patience time
\(1/\theta_i\), so its expected post-horizon holding and abandonment cost is
\(c_i/\theta_i+\gamma_i\).

For a horizon \(T>0\), the realized net value of a policy \(\pi\) is
\begin{equation}
  F^\pi(T)
  \coloneqq \sum_{m=1}^d r_mX_m^\pi(T)-\sum_{i=1}^n c_i\int_0^Tq_i^\pi(t)\,dt-\sum_{i=1}^n\gamma_iB_i^\pi(T)-C^{\mathrm{term}}\!\left(q^\pi(T)\right).
  \label{eq:obj}
\end{equation}
The decision maker seeks a policy that maximizes
\(\mathbb E[F^\pi(T)]\).

\section{Fluid Benchmark and Regret}
\label{sec:fluid-benchmark}
\subsection{Fluid Relaxation}
In this section we introduce a deterministic fluid relaxation of the stochastic
matching model and the corresponding linear optimization problem, whose optimal
value upper bounds \(\mathbb E[F^\pi(T)]/T\), the expected net value per unit
time achieved by any admissible policy \(\pi\in\Pi\).

In the fluid model, discrete arrivals and departures are replaced by continuous
flows.  Let
\(x(t)\in\mathbb{R}_+^{d}\) denote the vector of matching rates, and \(q(t)\in\mathbb{R}_+^{n}\) denote the vector of fluid queue
lengths, started at \(q(0)=0\).  Recalling that \(M\) is the incidence matrix
between agent types and match types, \((Mx(t))_i\) is the rate at which
type-\(i\) fluid is consumed by matching.

Type-\(i\) fluid enters at the constant rate \(\lambda_i\) and leaves through
two channels: it is matched at rate \((Mx(t))_i\), and it abandons at
intensity \(\theta_iq_i(t)\). The fluid queue
therefore evolves according to
\begin{equation}
    \dot q(t)=\lambda-Mx(t)-\Theta q(t),
    \qquad
    \Theta=\operatorname{diag}(\theta_1,\ldots,\theta_n),
    \label{eq:fluid-dynamics}
\end{equation}
and a matching-rate trajectory is admissible only if it keeps \(q(t)\ge0\)
throughout.

Setting \(\dot q(t)=0\) in \eqref{eq:fluid-dynamics}, a stationary point
\((x,q)\) satisfies the flow-balance condition
\begin{equation}
    Mx+\Theta q=\lambda.
    \label{eq:fluid-balance}
\end{equation}
This condition says that every agent who arrives eventually leaves, either
through a match or through abandonment.

We now define the optimization problem over the pairs \((x,q)\) satisfying
\eqref{eq:fluid-balance}.  Matching generates reward at rate \(r^\top x\), and
holding the queue costs \(c^\top q\) per unit time.  In addition, type-\(i\)
fluid abandons at rate \(\theta_iq_i\), and each abandoning unit incurs the
penalty \(\gamma_i\), so abandonment penalties are incurred at rate
\(\sum_{i=1}^n\gamma_i\theta_iq_i\).  The latter two costs are both linear in
\(q\), so we combine them into the abandonment-adjusted holding-cost vector
\begin{equation}
    \bar c_i=c_i+\theta_i\gamma_i,
    \qquad i\in[n],
    \label{eq:adjusted-holding-cost}
\end{equation}
which leaves \(r^\top x-\bar c^\top q\) as the objective, the net value
generated per unit time.  Maximizing it subject to \eqref{eq:fluid-balance}
yields the static fluid benchmark
\begin{align}
F^*(\lambda)=
\begin{bmatrix}
    \max_{x,q} \quad & \displaystyle\sum_{m=1}^d r_mx_m
    -\sum_{i=1}^n \bar c_iq_i \\
    \text{s.t.}\quad & \theta_iq_i+(Mx)_i=\lambda_i,
    \qquad i\in[n], \\
    & x\in\reals_{\geq 0}^d,\quad q\in\reals_{\geq 0}^n
\end{bmatrix}.
    \label{eq:fluidwithq}
\end{align}

The next lemma shows that the optimal value of \eqref{eq:fluidwithq} is an
upper bound on the net value per unit time achieved by any admissible policy
\(\pi\in\Pi\).

\begin{lemma}[Fluid Upper Bound]
\label{lmm:fluid-upper-bound}
For any admissible policy $\pi\in\Pi$,
\begin{equation*}
    \frac{\mathbb{E}[F^\pi(T)]}{T} \le F^*(\lambda).
\end{equation*}
\end{lemma}

The proof idea is as follows. Given \(\pi\in\Pi\), we use its
expected matching rates and queue lengths per unit time to construct a feasible
solution \((\bar{x},\bar{q})\) to \eqref{eq:fluidwithq}, whose objective
value equals \(\mathbb{E}[F^\pi(T)]/T\). The bound then follows from the
optimality of \(F^*(\lambda)\). The formal proof is deferred to Appendix~\ref{app:fluid-upper-bound}.

\subsection{High-Volume Setting}
\label{sec:high-volume-setting}
In this paper, we study the performance of matching policies under a high-volume
scaling in which the arrival rates of all agent types increase proportionally,
while all other model primitives remain fixed.  The scaling is expressed
relative to a base arrival-rate vector \(\lambda^b\in\mathbb R_+^n\) satisfying
\(\sum_{i=1}^n\lambda_i^b=1\).  We next define the regime formally.

\begin{definition}[High-volume regime]
\label{def:high-volume-regime}
We consider a family of systems indexed by
\(s\in\mathbb N\), where \(s\) tends to infinity, and let the arrival rates in
the \(s\)th system be
\begin{equation}
    \lambda_i^s=s\lambda_i^b
    \quad\text{for all }i\in\mathcal N,
    \qquad\text{or equivalently}\qquad
    \lambda^s=s\lambda^b.
    \label{eq:high-volume-scaling}
\end{equation}
Otherwise, each system \(s\) has the same basic structure as the system
described in Section~\ref{sec:model-setup}.
\end{definition}

In the high-volume regime, the fluid model in~\eqref{eq:fluid-dynamics} serves
as an accurate approximation to the discrete-event stochastic model of
Section~\ref{sec:model-setup}, whose dynamics are given by the flow balance
equation~\eqref{eq:flow-balance}. 
In addition to dynamic matching, similar fluid
approximations have also been studied for multiserver queues with abandonment;
see \cite{whitt2006fluid}.

\subsection{Performance Measure}
\label{sec:fluid-regret}

We evaluate matching policies by comparing their expected performance with
the fluid upper bound.  For each scale \(s\in\mathbb N\), let
\(F_s^\pi(T)\) denote the net value generated by policy \(\pi\) over the
horizon \([0,T]\) in the system with arrival-rate vector
\(\lambda^s=s\lambda^b\), and define
\begin{equation*}
    F_s^* \coloneqq F^*(\lambda^s).
\end{equation*}

Lemma~\ref{lmm:fluid-upper-bound} shows that the fluid value bounds from above
the expected net value per unit time of every admissible policy:
\begin{equation*}
    \frac{\mathbb E[F_s^\pi(T)]}{T}
    \leq F_s^*,
    \qquad
    s\in\mathbb N,\quad T>0,\quad \pi\in\Pi.
\end{equation*}
We therefore define regret as the gap between these two quantities.

\begin{definition}[Time-average regret]
\label{def:time-average-fluid-regret}
For a policy \(\pi\in\Pi\), its time-average regret at scale \(s\) over
horizon \(T\) is
\begin{equation*}
    \mathcal R_{T,s}^{\pi}
    \coloneqq
    F_s^*
    -
    \frac{\mathbb E[F_s^\pi(T)]}{T}
    \geq 0.
\end{equation*}
\end{definition}

We next introduce the performance measure used in this paper, the
approximation ratio, which compares a policy with the fluid benchmark.

\begin{definition}[Approximation ratio]
\label{def:approximation-ratio}
For a policy \(\pi\in\Pi\), its approximation ratio at scale \(s\) over
horizon \(T\) is
\begin{equation*}
    \rho_{T,s}^{\pi}
    \coloneqq
    \frac{\mathbb E[F_s^\pi(T)]}{TF_s^*}
    \leq 1,
\end{equation*}
where the fluid value \(F_s^*\) is assumed positive, and the inequality
follows from Lemma~\ref{lmm:fluid-upper-bound}.
\end{definition}

By Definition~\ref{def:time-average-fluid-regret},
\begin{equation*}
    \rho_{T,s}^{\pi}
    =
    1-\frac{\mathcal R_{T,s}^{\pi}}{F_s^*},
\end{equation*}
so a bound on the time-average regret \(\mathcal R_{T,s}^{\pi}\) directly
lower-bounds the approximation ratio.

\section{Main Results}
\label{sec:main-results}

In this section, we develop the longest-queue policy and state its
performance guarantee.  We first rewrite the fluid program~\eqref{eq:fluidwithq}
in a reduced form and derive its dual, which provides a convenient basis for the
policy design.  The \emph{general position gap} (GPG) condition, introduced next,
plays a central role in our main result.  Finally, we define
the longest-queue policy and state its performance guarantee.

Throughout the paper, we use \(\lambda\) for a general arrival-rate vector.
When a result depends explicitly on the high-volume scale \(s\), we write
\(\lambda^s=s\lambda^b\) and display that dependence explicitly.

\subsection{Reduced Formulation and Dual}
The fluid program~\eqref{eq:fluidwithq} can be written solely in terms of the
matching-rate vector \(x\), because the balance
constraint~\eqref{eq:fluid-balance} uniquely determines the queue-length vector
as \(q=\Theta^{-1}(\lambda-Mx)\), and the nonnegativity requirement \(q\geq0\)
is therefore equivalent to \(Mx\leq\lambda\).  Substituting this expression for
\(q\) into the objective of~\eqref{eq:fluidwithq} yields
\begin{align*}
    r^\top x-\bar c^\top q
    &=
    r^\top x-\bar c^\top\Theta^{-1}(\lambda-Mx) \notag\\
    &=
    -\lambda^\top\Theta^{-1}\bar c
    +\bigl(r+M^\top\Theta^{-1}\bar c\bigr)^\top x.
\end{align*}
For ease of
notation, we define the adjusted match-reward vector by
\begin{equation}
    \widetilde r
    =r+M^\top\Theta^{-1}\bar c,
    \qquad
    \widetilde r_m
    =r_m+\sum_{i=1}^n
    \frac{\bar c_i}{\theta_i}M_{im}.
    \label{eq:adjusted-match-reward}
\end{equation}
Using this notation, we obtain an alternative formulation of the linear
program~\eqref{eq:fluidwithq}:
\begin{align}
F^*(\lambda)=
\begin{bmatrix}
    \max_{x} \quad & \displaystyle\sum_{m=1}^d \widetilde r_mx_m
    -\sum_{i=1}^n\bar c_i\frac{\lambda_i}{\theta_i} \\
    \text{s.t.}\quad & (Mx)_i\leq\lambda_i,
    \qquad i\in[n], \\
    & x\in\reals_{\geq 0}^d
\end{bmatrix}.
    \label{eq:reduced-fluid-primal}
\end{align}
Then we derive the following dual formulation
of~\prettyref{eq:reduced-fluid-primal}:
\begin{align}
F^*(\lambda)=
\begin{bmatrix}
    \min_{U} \quad & \displaystyle U^\top\lambda
    -\sum_{i=1}^n\bar c_i\frac{\lambda_i}{\theta_i} \\
    \text{s.t.}\quad & \sum_{i=1}^n U_iM_{im}\geq\widetilde r_m,
    \qquad m\in[d], \\
    & U\in\reals_{\geq 0}^n
\end{bmatrix}.
    \label{eq:fluid-dual}
\end{align}

After eliminating \(q\), the matching-rate vector \(x\) is the only primal
decision variable, and \(U\) is the vector of dual multipliers associated
with the constraint \(Mx\leq\lambda\).  Let \(x^*\) be an optimal solution of
\eqref{eq:reduced-fluid-primal} and set \(q^*=\Theta^{-1}(\lambda-Mx^*)\);
then \((x^*,q^*)\) is optimal for \eqref{eq:fluidwithq}.  Throughout, we write
\(U^*\) for an optimal solution of the dual~\eqref{eq:fluid-dual}.
Since \(x^*\) is optimal for the reduced primal~\eqref{eq:reduced-fluid-primal},
the optimal objective value can be written as
\begin{equation}
    F^*(\lambda)
    =
    \langle \widetilde r,x^*\rangle
    -\sum_{i=1}^n\frac{\bar c_i\lambda_i}{\theta_i},
    \label{eq:fluid-value-at-x-star}
\end{equation}
which is also attained by \((x^*,q^*)\) in the original
formulation~\eqref{eq:fluidwithq}.

\subsection{General Position Gap}
Next, we introduce the notion of the \emph{general position gap} (GPG), which
measures how far the base arrival-rate vector \(\lambda^b\) can be perturbed
while the unique optimal solution of the dual~\eqref{eq:fluid-dual} remains
optimal.  We adopt the dual-based definition of \citet{wei2023constant}, which
differs slightly from the general position gap condition of
\citet{kerimov2023optimality}: the latter requires the primal
LP~\eqref{eq:reduced-fluid-primal} to have a unique nondegenerate optimal
solution, a stronger assumption.  Consequently, the condition of
\citet{kerimov2023optimality} implies a positive GPG defined
below, whereas the converse need not hold.

\begin{definition}[General Position Gap]\label{def:GPG}
Let \(\calB_\varepsilon(\lambda^b)\) denote the \(\ell_2\)-ball of radius
\(\varepsilon\) centered at \(\lambda^b\),
i.e., $\calB_\varepsilon(\lambda^b)
    =
    \{\widehat{\lambda}:
    \|\widehat{\lambda}-\lambda^b\|_2\leq\varepsilon\}$.
We say that \(\lambda^b\) has a GPG of at least \(\varepsilon>0\) if
\prettyref{eq:fluid-dual} has a unique optimal solution \(\dualopt\) and
\(\dualopt\) remains optimal when \(\lambda^b\) is replaced by any
\(\widehat{\lambda}\in\calB_\varepsilon(\lambda^b)\).
\end{definition}

A positive GPG ensures that 
the optimizer \(\dualopt\) of~\prettyref{eq:fluid-dual} is unique and remains
optimal throughout a neighborhood of \(\lambda^b\), so the match types it
selects stay correct even when the realized arrival rates fluctuate around
\(\lambda^b\).  We use this robustness to analyze the longest-queue policy proposed in \prettyref{sec:fg-policy}.

\subsection{The Longest-Queue Policy}
\label{sec:fg-policy}

We now describe our longest-queue policy, summarized in
Algorithm~\ref{alg:shadow-matching-neighbor-set}.
Before the system starts, we use the optimal dual solution \(U^*\)
of~\eqref{eq:fluid-dual} to identify the sets of \emph{good} matches
and \emph{critical} agent types, defined below.
These sets remain fixed throughout the horizon.
We say an agent type is \emph{critical} if its dual variable is positive,
and a match type is \emph{good} if its dual constraint is binding.
We define the sets of critical agent types and good match types as
\begin{equation}
    \mathcal I^+
    \coloneqq \{i\in\mathcal N:U_i^*>0\},
    \qquad
    \mathcal M^g
    \coloneqq
    \{m\in\mathcal M:\sum_{i=1}^n U_i^*M_{im}=\widetilde r_m\}.
    \label{eq:fg-critical-good}
\end{equation}
For a type-\(i\) arrival at time \(t\), define the sets of
available critical and noncritical good neighbors of \(i\) as
\begin{align}
    \mathcal N_i^+(t)
    &\coloneqq
    \{j\in\mathcal I^+:q_j(t)>0,\ (i,j)\in\mathcal M^g\},
    \label{eq:fg-available-critical-neighbors}\\
    \mathcal N_i^0(t)
    &\coloneqq
    \{j\in\mathcal N\setminus\mathcal I^+:
      q_j(t)>0,\ (i,j)\in\mathcal M^g\}\label{eq:fg-available-noncritical-neighbors}.
\end{align}
The system starts empty, \(q(0)=0\).
At each arrival time \(t\), let \(x(t)\) denote the matching decision.
For a type-\(i\) arrival, the policy selects
\(j^*\in\mathcal N_i^+(t)\) with the longest queue,
or \(j^*\in\mathcal N_i^0(t)\) if the former is empty,
and sets \(x(t)=e_{(i,j^*)}\).
If both sets are empty, it sets \(x(t)=0\), and the arrival
joins queue \(i\).\footnote{This ensures that our policy uses good matches whenever the agents are available, and thus no two waiting agents can ever form a good match. See~\prettyref{lmm:empty_neighbor} for the proof.}
Ties among longest queues and choices among noncritical types
follow a deterministic rule fixed in advance, such as choosing
the smallest index.
\begin{algorithm}[htbp]
\caption{Longest-Queue Policy}
\label{alg:shadow-matching-neighbor-set}
\algdef{SE}[EVENT]{Event}{EndEvent}[1]{\textbf{#1}}{}
\algtext*{EndEvent}
\algtext*{EndIf}
\begin{algorithmic}[1]
\Event{Upon a type-$i$ arrival at time $t$:}
    \If{$\mathcal N_i^+(t)\neq\varnothing$}
        \State Choose $j^*\in\operatorname*{arg\,max}_{j\in\mathcal N_i^+(t)}q_j(t)$ and $x(t)\gets e_{(i,j^*)}$. 
    \ElsIf{$\mathcal N_i^0(t)\neq\varnothing$}
        \State Choose $j^*\in\mathcal N_i^0(t)$ and $x(t)\gets e_{(i,j^*)}$. 
    \Else
        \State $x(t)\gets0$.
    \EndIf
    \State \refstepcounter{equation}\label{eq:fg-queue-update-arrival}$q(t+)=q(t)+e_i-Mx(t)$. \hfill(\theequation)
\EndEvent
\Event{Upon a type-$i$ abandonment at time $t$:}
    \State \refstepcounter{equation}\label{eq:fg-queue-update-abandonment}$q(t+)=q(t)-e_i$. \hfill(\theequation)
\EndEvent
\end{algorithmic}
\end{algorithm}

The intuition behind our policy is as follows. 
The optimal dual \(U_j^*\) measures the marginal (shadow) value of
type-\(j\) supply for the optimal adjusted matching reward.
For any match \(m=(i,j)\), the dual feasibility gives
\(\widetilde r_m-U_i^*-U_j^*\leq0\).
If the inequality is strict, the match earns less reward than
the shadow value it consumes, so the policy restricts attention
to good matches, for which equality holds.
Critical supply has positive shadow value, whereas
noncritical supply has zero shadow value.
The policy therefore prioritizes critical agents to limit
the loss of valuable supply through abandonment.
Among available critical partners, it serves the longest queue
to help keep critical queues small.

Next, we formally quantify how the policy's control of critical queues translates into performance relative to the fluid benchmark. Theorem~\ref{thm:average_regret} shows that, under the GPG condition, the longest-queue policy's time-average regret is bounded by a constant independent of the market size $s$ and operating horizon $T$, uniformly over $s\geq1$ and $T\geq1$. The bound also makes explicit its dependence on the number of agent types $n$ and the GPG parameter $\varepsilon$.

\begin{theorem}[Average Regret Bound]
\label{thm:average_regret}
Suppose that the base arrival-rate vector
\(\lambda^b\) satisfies the GPG condition in
Definition~\ref{def:GPG} with the parameter \(\varepsilon>0\).  Let
\(\pi=\pi^{\mathrm{LQ}}\) be the policy in
Algorithm~\ref{alg:shadow-matching-neighbor-set}.  Then, for every
\(s\geq1\) and \(T\geq1\),
\begin{equation*}
    0
    \;\leq\;
    \mathcal R_{T,s}^{\pi}
    \;\leq\;
    C\,\frac{n}{\varepsilon}\,\log\frac{n}{\varepsilon},
\end{equation*}
where \(C\) depends only on the abandonment rates \(\theta\) and the largest
adjusted match reward
\(\widetilde r_{\max}\) \(\bigl(\max_{m\in\mathcal M}\widetilde r_m\bigr)\).
Explicitly, for every \(s\geq1\) and \(T>0\),
\begin{equation*}
    \mathcal R_{T,s}^{\pi}
    =
    F_s^*
    -
    \frac{\mathbb E[F_s^\pi(T)]}{T}
    \leq
    \left(
        \frac{1}{T}
        + \theta_{\max}
    \right)
    \sqrt n\,\widetilde r_{\max}\,
    \mathcal B_q(s),
\end{equation*}
where $\theta_{\max}=\max_i \theta_i$ and
\begin{equation*}
\begin{aligned}
    \mathcal B_q(s)
   =\;&
    \frac{\sqrt{n}}{\varepsilon}
    +
    \left(
        1+
        \frac{4\sqrt{n}\left(\mathrm e-\frac{3}{2}\right)}
             {\varepsilon}
    \right)                                                     
    \times
    \ln\left(
        2
        +
        \frac{4\sqrt{n}}{\varepsilon}
        +
        \frac{
            4n 
        }{
            s\varepsilon^2
        } \sqrt{\sum_{i\in \calI^+} \theta^2_i}
    \right).
\end{aligned}
\end{equation*}
In particular, \(\mathcal B_q(s)\) is non-increasing in \(s\), and
\begin{equation*}
    \sup_{s\geq1}\mathcal B_q(s)
    =
    \mathcal B_q(1)
    =
    O\left(
        \frac{\sqrt{n}}{\varepsilon}
        \log\frac{\sqrt{n}}{\varepsilon}
    \right).
\end{equation*}
\end{theorem}
The proof of Theorem~\ref{thm:average_regret} is 
presented in Section \ref{sec:analysis-lq}.
Note that Theorem~\ref{thm:average_regret} immediately leads to the following lower bound on the approximation ratio.

\begin{corollary}[Approximation Ratio]
\label{cor:approx_ratio}
Suppose the assumptions of Theorem~\ref{thm:average_regret} hold.
Then
\begin{equation*}
    \rho_{T,s}^{\pi^{\mathrm{LQ}}}
    \geq
    1- \frac{C}{s},
\end{equation*}
where 
$C$ is a constant independent of $T$ and $s$. It follows that the approximation ratio approaches $1$ at rate \(1/s\) as $s \to \infty$. 
\end{corollary}

Without the GPG condition, our longest-queue policy still achieves an
approximation ratio that approaches \(1\), with a convergence guarantee
of order \(1/\sqrt{s}\) rather than \(1/s\). We present the formal result
in Proposition~\ref{prop:no_gpg}.

\begin{proposition}[Regret and Approximation Ratio without GPG]
\label{prop:no_gpg}
For any optimal solution \(U^*\) of
\eqref{eq:fluid-dual} at the base arrival-rate vector \(\lambda^b\), let
\(\pi=\pi^{\mathrm{LQ}}\) be the policy in
Algorithm~\ref{alg:shadow-matching-neighbor-set}. Then, for every
\(s\geq1\) and \(T>0\),
\begin{equation}
    0\leq\mathcal R_{T,s}^{\pi}
    \leq
    \left(\frac{1}{T}+\|\theta\|_\infty\right)
    \sqrt n\,\widetilde r_{\max}\sqrt{\frac{s}{\theta_{\min}}},
    \label{eq:no-gpg-regret}
\end{equation}
where \(\theta_{\min}=\min_{i\in\mathcal N}\theta_i\) and
\(\widetilde r_{\max} =\max_{m\in\mathcal M}\widetilde r_m\).
Moreover, for every \(s\geq1\) and \(T\geq1\),
\begin{equation}
    \rho_{T,s}^{\pi}
    \geq 1-\frac{C}{\sqrt{s}},
    \label{eq:no-gpg-ratio}
\end{equation}
where \(C\) is a constant independent of \(T\) and \(s\).
\end{proposition}

The proof is given in Appendix~\ref{app:proof_no_gpg}. 
Without GPG, matching contributes only a
nonpositive drift term, while abandonment supplies the downward drift that
yields the \(O(\sqrt{s})\) bound on
\(\mathbb E[\|q(t)\|_{\mathcal I^+}]\).
This explains why the regret bound in Proposition~\ref{prop:no_gpg} depends on \(1/\sqrt{\theta_{\min}}\) and
does not extend directly to the non-abandonment setting.

Under GPG, the longest-queue rule
itself generates a stronger matching-induced drift, improving the queue bound
to \(O(1)\) and the approximation gap from \(O(1/\sqrt{s})\) to \(O(1/s)\). Thus,  when there is no abandonment (i.e., $\theta_{\max}=0$), Theorem \ref{thm:average_regret}
yields \(O(1/T)\) time-average regret, or equivalently, \(O(1)\) cumulative
regret. 

For comparison, \citet{aveklouris2021matching} develop matching-rate-based
and priority-ordering policies in a similar high-volume regime. Both policies make matching decisions only at scheduled review times,
with intervals proportional to \(s^{-2/3}\). Let
\(\pi^{\mathrm r}\) denote either policy. Expressed in our approximation-ratio
notation, their results imply
\[
    \lim_{T\to\infty}\lim_{s\to\infty}
    \rho_{T,s}^{\pi^{\mathrm r}}=1.
\]
Their proof first takes \(s\to\infty\), so that the large number
of arrivals within each review period yields a fluid approximation
of the system.
 Then \(T\to\infty\) is taken so that the fluid queues approach
the invariant queue lengths specified by their optimal fluid
benchmark. 

In contrast,
Corollary~\ref{cor:approx_ratio} and Proposition~\ref{prop:no_gpg} provide
explicit approximation guarantees at every market size \(s\geq1\) and
horizon \(T\geq1\), without taking either limit.

\section{Analysis of the Longest-Queue Policy}
\label{sec:analysis-lq}

In this section, we prove Theorem~\ref{thm:average_regret}
for the longest-queue policy defined in
Algorithm~\ref{alg:shadow-matching-neighbor-set}.
We first introduce a restricted inner product and
norms that will be used throughout the rest of the paper.  For vectors \(u,v\in\mathbb R^n\), write
\[
    \langle u,v\rangle_{\mathcal I^+}
    \coloneqq
    \sum_{i\in\mathcal I^+}u_iv_i
\]
and
\[
    \|u\|_{\mathcal I^+}
    \coloneqq
    \left(\sum_{i\in\mathcal I^+} u_i^2\right)^{1/2},
    \qquad
    \|u\|_{1,\mathcal I^+}
    \coloneqq
    \sum_{i\in\mathcal I^+}|u_i|
\]
for the corresponding \(\ell_2\)- and \(\ell_1\)-norms. Unless otherwise specified, \(\|\cdot\|\) denotes the Euclidean
\(\ell_2\)-norm, \(\|\cdot\|_1\) denotes the \(\ell_1\)-norm, and
\(\| u \|_\infty \coloneqq \max_{i}| u_i|\) the maximum norm.

The proof proceeds in four steps.  First, in \prettyref{sec:step1}, we express
the regret \(TF^*-\mathbb E[F^\pi(T)]\) in terms of the queue lengths of the
critical agents, \(q_i(t)\) for \(i\in\mathcal I^+\). Second, in \prettyref{sec:step2}, we establish a
negative-drift condition for \(\cnorm{q(t)}^2\).
Third, in \prettyref{sec:step3}, we use this drift
condition to obtain a bound on
\(\mathbb E[\cnorm{q(t)}]\) that is uniform in
the time \(t\) and in the scale \(s\). Finally, in \prettyref{sec:step4}, we
substitute this queue bound into the regret representation to prove
Theorem~\ref{thm:average_regret}.

\subsection{Bounding Regret Using Critical Queue Lengths}\label{sec:step1}
In this subsection, we express regret in terms of the critical queue lengths. First, although the objective nets matching rewards against
holding and abandonment costs, Lemma~\ref{lmm:bound_reward} shows that the
regret $T \cdot F^* - \E[F^{\pi}(T)]$ reduces to a gap in adjusted reward.
Specifically, it equals the difference between the fluid benchmark's
cumulative adjusted reward $T\langle \tilde{r}, x^* \rangle$ and the policy's
expected cumulative adjusted reward
 $\E[\sum_{m} \tilde{r}_m X_m^\pi(T)]$.  

\begin{lemma}\label{lmm:bound_reward}
For any admissible policy $\pi\in\Pi$, the regret satisfies
\begin{align*}
    T F^*-\mathbb E[F^\pi(T)]
    =
    T\langle\tilde r,x^*\rangle
    -
    \mathbb E\!\left[
        \sum_{m=1}^d\tilde r_m X_m^\pi(T)
    \right],
\end{align*}
where the adjusted reward $\tilde r_m$ is defined
in~\eqref{eq:adjusted-match-reward}.
\end{lemma}

The proof of Lemma~\ref{lmm:bound_reward} is deferred to Appendix~\ref{app:lmm:bound_reward}.
For the remainder of the analysis, we focus on our longest-queue policy
\(\pi=\pi^{\mathrm{LQ}}\), so we suppress the superscript \(\pi\) on the
processes \(X^\pi(t)\), \(B^\pi(t)\), and \(q^\pi(t)\) for notational
simplicity. With a slight abuse of notation, \(\theta q\) denotes the
element-wise product \((\theta_1q_1,\ldots,\theta_nq_n)\), and similarly
for \(\theta q(t)\) and \(\theta q^*\).

Next, we show in Proposition~\ref{prop:regret_decomposition} that
the regret can be further expressed entirely in terms of the critical queue
lengths \(q_i(t)\), \(i\in\mathcal I^+\), which follows from complementary
slackness and the policy's restriction to good matches.

\begin{proposition}
\label{prop:regret_decomposition}
Under the longest-queue policy $\pi$, the regret satisfies
\begin{align}
    T F^* - \E[F^\pi(T)] 
    &= \langle U^*, \E[q(T)] \rangle_{\mathcal{I}^+} + \int_0^T \E\left[\langle U^*, \theta q(t) \rangle_{\mathcal{I}^+} \right] dt. \label{eq:regret_decomp}
\end{align}
\end{proposition}

\begin{proof}
By Lemma~\ref{lmm:bound_reward}, the regret satisfies
 \begin{align}
    T F^* - \E[F^\pi(T)] &= \langle \tilde{r}, T x^* - \E\left[ X(T)\right] \rangle \notag \\
    &= \langle M^\top U^*, T x^* - \E\left[ X(T)\right] \rangle
    = \langle  U^*, T Mx^* - \E\left[ M X(T)\right] \rangle,
    \label{eq:regret-dual-form}
\end{align}
where the second equality holds as follows. For $m\notin\mathcal M^g$, we have
$(M^\top U^*)_m > \tilde{r}_m$. Complementary slackness
for~\eqref{eq:reduced-fluid-primal} and~\eqref{eq:fluid-dual} gives
$x^*_m=0$, and our policy restricts match selection to the good matches
$\mathcal M^g$, so $X_m(T)=0$, hence the $m$-th component of
$T x^* - \E[X(T)]$ vanishes. For $m\in\mathcal M^g$,
$(M^\top U^*)_m=\tilde{r}_m$ by the definition of $\mathcal M^g$
in~\eqref{eq:fg-critical-good}. Replacing $\tilde{r}$ by $M^\top U^*$ thus
leaves the inner product unchanged.

Now, taking expectations in the flow-balance
equation~\eqref{eq:flow-balance}
and using the abandonment
identity
(see~\eqref{eq:expected-abandonments} in the proof of
Lemma~\ref{lmm:fluid-upper-bound}
in Appendix~\ref{app:fluid-upper-bound})
$$
 \mathbb{E}[B(T)]
    =
    \int_0^T\mathbb{E}[\theta q(t)]\,dt,
$$    
together with $\E[A(T)]=\lambda T$ for the
Poisson arrivals, we have
\begin{equation}
    \E[MX(T)]
    =
    \E[A(T)-B(T)-q(T)]
    =
    \lambda T
    - \int_0^T \E[\theta q(t)]\,dt
    - \E[q(T)].
    \label{eq:expected-matched-mass}
\end{equation}
By the equality constraint in the fluid model~\eqref{eq:fluidwithq}, the
optimal solution satisfies $Mx^* = \lambda - \theta q^*$, so
$T Mx^* = \lambda T - \theta q^* T$.  Substituting this
and~\eqref{eq:expected-matched-mass} into~\eqref{eq:regret-dual-form}, the
$\lambda T$ terms cancel and we obtain
\begin{align*}
  T F^* - \E[F^\pi(T)]
  & =
  \iprod{U^*}{\int_0^T \expect{   \theta q(t)} dt - \theta q^* T + \expect{q(T)}} \\
 &  = \iprod{U^*}{\int_0^T \expect{   \theta q(t)} dt - \theta q^* T + \expect{q(T)}}_{\calI^+} \\
& = \iprod{U^*}{\int_0^T  \expect{   \theta q(t)} dt}_{\calI^+} 
+ \iprod{U^*}{\expect{q(T)}}_{\calI^+},
\end{align*}
where the second equality holds because $U_i^*=0$ for $i\notin\mathcal I^+$
by the definition of $\mathcal I^+$ in~\eqref{eq:fg-critical-good}, so only
the coordinates in $\mathcal I^+$ contribute; the last equality holds
because $U_i^*>0$ for $i\in\mathcal I^+$, so complementary slackness
for~\eqref{eq:reduced-fluid-primal} and~\eqref{eq:fluid-dual} forces
$(Mx^*)_i=\lambda_i$, and the fluid balance~\eqref{eq:fluid-balance} then
gives $\theta_i q_i^* = \lambda_i - (Mx^*)_i = 0$.
\end{proof}

\subsection{Negative Drift Analysis of $\norm{q(t)}_{\mathcal{I}^+}^2$}\label{sec:step2}

Since the regret can be expressed in terms of the critical queue
lengths by~\eqref{eq:regret_decomp} of
Proposition~\ref{prop:regret_decomposition}, bounding regret
reduces to bounding $\mathbb{E}[\|q(t)\|_{\mathcal I^+}]$.
To obtain such a bound, we first establish a negative drift
condition for $\|q(t)\|_{\mathcal I^+}^2$.
We begin by defining the infinitesimal generator, which characterizes
the instantaneous evolution of a continuous-time Markov chain.

\begin{definition}[Infinitesimal Generator]
\label{def:generator}
Let $\mathbf{Y}$ be a state space and let $\{Y(t)\}_{t \ge 0}$ be a time-homogeneous continuous-time Markov chain on $\mathbf{Y}$. Following \cite[Section~2.4.1]{bujorianu2012stochastic}, the \emph{infinitesimal generator} $\mathcal{L}$ acting on a bounded measurable function $f$ is defined by
\begin{align*}
    \mathcal{L} f(y) := \lim_{t \downarrow 0} \frac{\mathbb{E}[f(Y(t)) \mid Y(0) = y] - f(y)}{t}, \quad y \in \mathbf{Y},
\end{align*}
whenever this limit exists.
\end{definition}

We now state the main result of this subsection, which establishes
a negative drift condition for $\|q(t)\|_{\mathcal I^+}^2$
under the GPG condition.
To display the dependence on the scale $s$ explicitly, we write
$\lambda^s=s\lambda^b$ in place of the general arrival-rate
vector $\lambda$ throughout the following analysis.
\begin{proposition}\label{prop:queue_square_negdrift}
    Assume that $\lambda^b$ satisfies the GPG condition with $\epsilon$. Under the longest-queue policy described in Algorithm~\ref{alg:shadow-matching-neighbor-set}, given $q(t)$, we have
    \begin{align*}
 \mathcal{L}\|q(t)\|^2_{\mathcal{I}^+} \leq  2s \left( -\frac{\epsilon}{\sqrt{n}} \|q(t)\|_{\mathcal{I}^+} \right) + s  . \notag
\end{align*}
\end{proposition}

We prove Proposition~\ref{prop:queue_square_negdrift}
by adapting the Lyapunov-drift argument of
\citet{kerimov2023optimality}, with three modifications. First, their model has no abandonment and evolves in discrete time, whereas our drift is computed through the infinitesimal generator of a continuous-time process; abandonment contributes an additional state-dependent term, which we show is nonpositive. Second, our critical types \(\mathcal I^+\) and good matches \(\mathcal M^g\) are defined through an optimal dual solution \(U^*\), whereas their corresponding sets are identified from a unique nondegenerate primal solution. Third, they remove noncritical agents that are not matched immediately, so the corresponding queues are identically zero. We retain these agents and must account for their presence when analyzing the drift of the critical queues. 

We begin by introducing notation for the matching decisions.
At an arrival epoch, the policy selects a decision
$x(t)\in\{0,1\}^d$ with at most one nonzero entry.
Conditional on a type-$i$ arrival, we write $x^i(t)$ for this
decision, which is deterministic given $q(t)$.
We define the aggregate matching rates
\begin{equation}
    \bar{x}(t) \coloneqq \sum_{i=1}^n \lambda_i^s x^i(t),
    \qquad
    \bar{x}^b(t) \coloneqq \sum_{i=1}^n \lambda_i^b x^i(t),
    \label{eq:aggregate-matching-rate}
\end{equation}
so that $\bar{x}(t)=s\bar{x}^b(t)$ under the high-volume
scaling~\eqref{eq:high-volume-scaling}.

As the first step in proving
Proposition~\ref{prop:queue_square_negdrift}, we calculate the
infinitesimal generator acting on $\|q(t)\|_{\mathcal I^+}^2$
and obtain the following bound
in terms of the aggregate matching rate $\bar{x}^b(t)$.

\begin{lemma}\label{lmm:generator_bound}
Under the longest-queue policy described in Algorithm~\ref{alg:shadow-matching-neighbor-set}, given $q(t)$, we have
\begin{align*}
    \mathcal{L}\|q(t)\|^2_{\mathcal{I}^+} \leq 2s\langle q(t), \lambda^b - M\bar{x}^b(t) \rangle_{\mathcal{I}^+} + s.
\end{align*}
\end{lemma}

The proof of Lemma~\ref{lmm:generator_bound} is deferred to
Appendix~\ref{app:lmm:generator_bound}.

The bound in Lemma~\ref{lmm:generator_bound} reduces the proof of
Proposition~\ref{prop:queue_square_negdrift} to showing that
\begin{align}
    \langle q(t),\lambda^b-M\bar{x}^b(t)\rangle_{\mathcal I^+}
    \leq -\frac{\varepsilon}{\sqrt n}
    \|q(t)\|_{\mathcal I^+}.
    \label{eq:drift_ineq}
\end{align}

We introduce some notation. With a slight abuse of notation, we write $(i, j)$ for the match $m$ between types $i$ and $j$. Let
\begin{align*}
    \mathcal{Q}^+(t) := \{i \in \mathcal{N} : q_i(t) > 0\}
\end{align*}
denote the set of agent types with positive queue length at time $t$. Define
\begin{align*}
    \mathcal{M}^+(t) := \{ m = (i,j) \in \mathcal{M}^g : \{i, j\} \cap \mathcal{Q}^+(t) \neq \emptyset \}
\end{align*}
as the set of good matches involving at least one type with a positive queue.

We prove inequality \eqref{eq:drift_ineq} in two steps using an auxiliary variable $\xi \in \mathcal{X}(t)$, where 
    \begin{equation}\label{eq:auxiliary-feasible-set}
    \mathcal X(t)
    := \left\{x\in\mathbb R_+^d:
    \begin{array}{l}
        x_m=0 \text{ for all } m\notin\mathcal M^+(t)
        \text{ and} \\[0.5ex]
        (Mx)_i\leq\lambda_i^b
        \text{ for all } i\in[n]\setminus\mathcal Q^+(t)
    \end{array}
    \right\}.
\end{equation}
First, we establish that $\bar{x}^b(t)$ maximizes
$\langle q(t),Mx\rangle_{\mathcal I^+}$ over the feasible set
$\mathcal X(t)$, which implies that
\begin{align*}
    \langle q(t),\lambda^b-M\bar{x}^b(t)\rangle_{\mathcal I^+}
    \leq
    \langle q(t),\lambda^b-M\xi\rangle_{\mathcal I^+}
\end{align*}
for any feasible auxiliary vector $\xi\in\mathcal X(t)$.
This optimality result is established in
Lemma~\ref{lmm:optimal_action}.
\begin{lemma}\label{lmm:optimal_action}
Under the longest-queue policy described in
Algorithm~\ref{alg:shadow-matching-neighbor-set},
$\bar{x}^b(t)$ is an optimal solution of the following problem:
\begin{align}\label{eq:max_queue_prop}
    \max_{x\in\mathcal X(t)}
    \langle q(t),Mx\rangle_{\mathcal I^+},
\end{align}
where $\mathcal X(t)$ is defined
in~\eqref{eq:auxiliary-feasible-set}.
\end{lemma}

The proof of Lemma~\ref{lmm:optimal_action} is deferred to
Appendix~\ref{app:lmm:optimal_action}. Second, we use the GPG assumption to construct a feasible
auxiliary vector $\xi\in\mathcal X(t)$ satisfying
\begin{align*}
    \langle q(t),\lambda^b-M\xi\rangle_{\mathcal I^+}
    \leq -\frac{\varepsilon}{\sqrt n}
    \|q(t)\|_{\mathcal I^+}.
\end{align*}
The existence of such a vector is established in
Lemma~\ref{lmm:construct_aux}.

\begin{lemma}\label{lmm:construct_aux}
Assume that $\lambda^b$ satisfies the GPG condition with $\epsilon>0$. Then there exists
$\xi\in\mathcal{X}(t)$ such that
\begin{align}
    \langle q(t),\lambda^b-M\xi\rangle_{\mathcal{I}^+}
    \leq -\frac{\epsilon}{\sqrt n}
    \|q(t)\|_{\mathcal{I}^+}.
    \label{eq:drift_construct}
\end{align}
\end{lemma}

The proof of Lemma~\ref{lmm:construct_aux} is deferred to
Appendix~\ref{app:lmm:construct_aux}.
We now combine the preceding bounds to complete the proof of
Proposition~\ref{prop:queue_square_negdrift}.

\begin{proof}[Proof of Proposition~\ref{prop:queue_square_negdrift}]
Choose $\xi\in\mathcal X(t)$ as in
Lemma~\ref{lmm:construct_aux}.
By Lemmas~\ref{lmm:generator_bound} and~\ref{lmm:optimal_action},
\begin{align*}
    \mathcal L\|q(t)\|_{\mathcal I^+}^2
    &\leq
    2s\langle q(t),\lambda^b-M\bar{x}^b(t)\rangle_{\mathcal I^+}+s\\
    &\leq
    2s\langle q(t),\lambda^b-M\xi\rangle_{\mathcal I^+}+s\\
    &\leq
    -\frac{2s\epsilon}{\sqrt n}\|q(t)\|_{\mathcal I^+}+s,
\end{align*}
where the last inequality follows from~\eqref{eq:drift_construct}.
\end{proof}

\subsection{Bound on $\cnorm{q(t)}$}\label{sec:step3}

Using the negative drift condition established in
Section~\ref{sec:step2}, we derive a uniform bound on the
expected critical queue length
$\mathbb{E}[\|q(t)\|_{\mathcal I^+}]$ for all $t \geq 0$,
as stated in the following proposition.
\begin{proposition}\label{prop:queue_length_bound}
    Assume that $\lambda^b$ satisfies the GPG condition with parameter $\varepsilon > 0$. Under Algorithm~\ref{alg:shadow-matching-neighbor-set}, the expected queue length satisfies the uniform bound
    \begin{align*}
    \mathbb{E}\big[\|q(t)\|_{\mathcal{I}^+}\big] \le \frac{\sqrt{n}}{\varepsilon} + \left( 1 + \frac{2\sqrt{n}(2e - 3)}{\varepsilon} \right) \ln \left( 2 + \frac{4\sqrt{n}}{\varepsilon} + \frac{4n\|\theta\|_{\mathcal{I}^+}}{s\varepsilon^2} \right).
    \end{align*}
\end{proposition}

To prove Proposition~\ref{prop:queue_length_bound}, we establish two properties of $\|q(t)\|_{\mathcal I^+}$.
The first is that it has bounded jumps, as stated in the
following lemma.

\begin{lemma}\label{lmm:bounded_jump}
Under Algorithm~\ref{alg:shadow-matching-neighbor-set}, the one-step change in the critical queue norm at any event time $t$ satisfies
\begin{align*}
\left| \|q(t+)\|_{\mathcal{I}^+} - \|q(t)\|_{\mathcal{I}^+} \right| \leq 1.
\end{align*}
\end{lemma}

The proof of Lemma~\ref{lmm:bounded_jump} follows from a
direct check of the change in $\|q(t)\|_{\mathcal I^+}$
at each arrival or abandonment. The details are provided
in Appendix~\ref{app:lmm:bounded_jump}.

The second property is that $\|q(t)\|_{\mathcal I^+}$ has
negative drift when it is sufficiently large, as stated
in the following lemma.

\begin{lemma}\label{lmm:queue_negdrift}
Assume that $\lambda^b$ satisfies the GPG condition with parameter $\varepsilon > 0$. Under Algorithm~\ref{alg:shadow-matching-neighbor-set}, the generator of $\|q(t)\|_{\mathcal{I}^+}$ satisfies
\begin{align*}
\mathcal{L}\|q(t)\|_{\mathcal{I}^+} \leq 
\begin{cases}
-\dfrac{s\varepsilon}{2\sqrt n} & \text{if } \|q(t)\|_{\mathcal{I}^+} \geq \dfrac{\sqrt{n}}{\varepsilon}, \\[10pt]
s + \dfrac{\sqrt{n}\|\theta\|_{\mathcal{I}^+}}{\varepsilon} & \text{if } \|q(t)\|_{\mathcal{I}^+} < \dfrac{\sqrt{n}}{\varepsilon}.
\end{cases}
\end{align*}
\end{lemma}

The proof of Lemma~\ref{lmm:queue_negdrift} uses the
squared-norm drift bound in
Proposition~\ref{prop:queue_square_negdrift}, with the
full derivation given in Appendix~\ref{app:lmm:queue_negdrift}.

Having established bounded jumps and a negative drift
condition for $\|q(t)\|_{\mathcal I^+}$ in
Lemmas~\ref{lmm:bounded_jump} and~\ref{lmm:queue_negdrift},
respectively, we next derive a uniform bound on
$\mathbb E[\|q(t)\|_{\mathcal I^+}]$.
To obtain this bound, we first establish a general result for
continuous-time Markov jump processes in
Lemma~\ref{lmm:general_drift} below.

\begin{lemma}\label{lmm:general_drift}
Let $\{Y(t)\}_{t \geq 0}$ be a continuous-time Markov process on a state space $\mathcal{Y}$ with jump types $\mathcal{J}$, where $\lambda_j(y)$ denotes the rate of jump $j \in \mathcal{J}$ from state $y$ and $\Delta_j(y)$ the associated change in the state, so that a jump of type $j$ moves the process from $y$ to $y + \Delta_j(y)$. For a real-valued function $f$ on $\mathcal{Y}$, write $\Delta_j f(y) := f(y + \Delta_j(y)) - f(y)$ for the change in $f$ produced by such a jump.
Assume the following conditions hold:
\begin{enumerate}
    \item \label{cond:jump-rep} Jump representation: the generator of Definition~\ref{def:generator} satisfies
    \begin{align*}
        \mathcal{L}\varphi(y) = \sum_{j \in \mathcal{J}} \lambda_j(y) \Delta_j \varphi(y)
    \end{align*}
    for $\varphi = f$ and for $\varphi = e^{\alpha f}$ with $\alpha > 0$.
    \item \label{cond:drift} Drift condition: There exist constants $\delta > 0$, $D \geq 0$, and $B \geq 0$ such that
    \begin{align*}
        \mathcal{L}f(y) \leq -\delta \quad \text{if } f(y) \geq D, \quad \text{and} \quad \mathcal{L}f(y) \leq B \quad \text{if } f(y) < D.
    \end{align*}
    \item \label{cond:bounded-jump} Bounded jump size: $|\Delta_j f(y)| \leq K$ for all $j \in \mathcal{J}$ and $y \in \mathcal{Y}$.
    \item \label{cond:bounded-rate} Bounded upward rate: The aggregate rate of positive jumps satisfies
    \begin{align*}
        \sum_{j \in \mathcal{J}^+(y)} \lambda_j(y) \leq \Lambda,
    \end{align*}
    where $\mathcal{J}^+(y) = \{j \in \mathcal{J} : \Delta_j f(y) > 0\}$.
\end{enumerate}
Then for any $t \geq 0$, we have
\begin{align}
    \mathbb{E}[f(Y_t)] \leq \frac{1}{\alpha^*} \max\left\{ \alpha^* f(Y_0),\, \ln \frac{B + \delta}{\delta - \alpha^* C_1} + \alpha^* D \right\}, \label{eq:general_drift_bound}
\end{align}
where $C_1 = \frac{K\delta}{2} + (e - \frac{3}{2})\Lambda K^2 $, $\alpha^* = \min\left( \frac{1}{K}, \frac{\delta}{2C_1} \right)$ and $Y_0=Y(0)$ is the deterministic initial state
of the process.
\end{lemma}

We outline the proof idea of
Lemma~\ref{lmm:general_drift} below, with the full proof
provided in Appendix~\ref{app:lmm:general_drift}.
Following \citet{hajek1982hitting}, our proof uses quadratic bounds on exponential increments to convert negative drift of a function \(f\) into negative drift of \(e^{\alpha f}\). We adapt this argument to the continuous-time generator while allowing the aggregate downward jump rate to be unbounded. For downward jumps, bounded jump size permits the quadratic remainder to be absorbed into the negative first-order contribution, so no bound on their aggregate rate is needed. For upward jumps, the quadratic remainder is controlled by assuming that their aggregate rate is bounded. Choosing \(\alpha>0\) sufficiently small then yields negative drift of the exponential Lyapunov function outside a bounded region and, consequently, a uniform exponential-moment bound. Jensen’s inequality then gives the desired uniform bound on \(\mathbb E[\|q(t)\|_{\mathcal I^+}]\).

We now apply Lemma~\ref{lmm:general_drift} to
$\|q(t)\|_{\mathcal I^+}$ to prove
Proposition~\ref{prop:queue_length_bound}.

\begin{proof}[Proof of Proposition~\ref{prop:queue_length_bound}]
We apply Lemma~\ref{lmm:general_drift} with $Y(t) = q(t)$ and $f(q) = \|q\|_{\mathcal{I}^+}$. The jump types are arrivals (type $i$, rate $\lambda_i^s$) and abandonments (type $i \in \mathcal{I}^+$, rate $\theta_i q_i$), with corresponding changes
\begin{align*}
    \Delta_{\mathrm{arr},i} f(q) = \|q + e_i - Mx^i\|_{\mathcal{I}^+} - \|q\|_{\mathcal{I}^+}, \qquad \Delta_{\mathrm{abd},i} f(q) = \|q - e_i\|_{\mathcal{I}^+} - \|q\|_{\mathcal{I}^+}.
\end{align*}

First, we verify the conditions of Lemma~\ref{lmm:general_drift}.
Condition~\ref{cond:jump-rep} is the generator formula for
the queue process $q(t)$, applied to $f$ and $e^{\alpha f}$
(see Lemma~\ref{lmm:generator}).
Condition~\ref{cond:drift} holds by Lemma~\ref{lmm:queue_negdrift}, with
\begin{align*}
    \delta = \frac{s\varepsilon}{2\sqrt{n}}, \qquad D = \frac{\sqrt{n}}{\varepsilon}, \qquad B = s + \frac{\sqrt{n}\|\theta\|_{\mathcal{I}^+}}{\varepsilon}.
\end{align*}
Condition~\ref{cond:bounded-jump} holds with $K = 1$ by Lemma~\ref{lmm:bounded_jump}, and Condition~\ref{cond:bounded-rate} holds with $\Lambda = s$, since only arrivals can increase $\|q\|_{\mathcal{I}^+}$ and the total arrival rate is $\sum_{i=1}^n \lambda_i^s = s$.

Next, we compute the constant $C_1$. Substituting the values above, we obtain
\begin{align*}
    C_1 = \frac{K\delta}{2} + (e - \frac{3}{2})\Lambda K^2  = \frac{s\varepsilon}{4\sqrt{n}} + s(e - \frac{3}{2}).
\end{align*}
To determine $\alpha^* = \min\left( \frac{1}{K}, \frac{\delta}{2C_1} \right)$, we compute
\begin{align*}
    \frac{\delta}{2C_1} = \frac{\frac{s\varepsilon}{2\sqrt{n}}}{\frac{s\varepsilon}{2\sqrt{n}} + 2s(e - \frac{3}{2})} < 1 = \frac{1}{K},
\end{align*}
so $\alpha^* = \frac{\delta}{2C_1}$. The inverse of $\alpha^*$ simplifies to
\begin{align*}
    \frac{1}{\alpha^*} = \frac{2C_1}{\delta} = 1 + \frac{2\sqrt{n}(2e - 3)}{\varepsilon}.
\end{align*}
Since $\alpha^* C_1 = \frac{\delta}{2}$, the logarithmic term in~\eqref{eq:general_drift_bound} simplifies to
\begin{align*}
    \frac{B + \delta}{\delta - \alpha^* C_1} = \frac{B + \delta}{\delta/2} = 2\left(\frac{B}{\delta} + 1\right) = 2 + \frac{4\sqrt{n}}{\varepsilon} + \frac{4n\|\theta\|_{\mathcal{I}^+}}{s\varepsilon^2}.
\end{align*}
Since $q(0) = 0$, the bound~\eqref{eq:general_drift_bound} from Lemma~\ref{lmm:general_drift} yields
\begin{align*}
    \mathbb{E}\bigl[\|q(t)\|_{\mathcal{I}^+}\bigr] \leq \frac{\sqrt{n}}{\varepsilon} + \left( 1 + \frac{2\sqrt{n}(2e - 3)}{\varepsilon} \right) \ln \left( 2 + \frac{4\sqrt{n}}{\varepsilon} + \frac{4n\|\theta\|_{\mathcal{I}^+}}{s\varepsilon^2} \right). 
\end{align*}
\end{proof}

\subsection{Proof of the Main Result}\label{sec:step4}

We now combine the regret representation
in Proposition~\ref{prop:regret_decomposition} and the uniform
queue bound in Proposition~\ref{prop:queue_length_bound}
to prove Theorem~\ref{thm:average_regret}.
The GPG condition enters through the negative-drift argument
underlying Proposition~\ref{prop:queue_length_bound}.

\begin{proof}[Proof of Theorem~\ref{thm:average_regret}]
    Consider the system with arrival-rate vector
\(\lambda^s=s\lambda^b\), writing
    \(F_s^*\), \(F_s^\pi(T)\), and \(\mathcal{B}_q(s)\) for the corresponding
    quantities.
    The nonnegativity of the regret follows from the fluid upper bound in
    Lemma~\ref{lmm:fluid-upper-bound}.  For the upper bound, by
    Proposition~\ref{prop:regret_decomposition}, we have
    \begin{align*}
        T F_s^* - \E[F_s^\pi(T)] 
        = \langle U^*, \E[q(T)] \rangle_{\mathcal{I}^+} + \int_0^T \E\left[\langle U^*, \theta q(t) \rangle_{\mathcal{I}^+} \right] \diff t.
    \end{align*}
    Applying the Cauchy--Schwarz inequality to each term yields
    \begin{align*}
        T F_s^* - \E[F_s^\pi(T)]   
        &\leq \|U^*\|_{\mathcal{I}^+} \E\bigl[\|q(T)\|_{\mathcal{I}^+}\bigr] + \int_0^T \|\theta\|_\infty \|U^*\|_{\mathcal{I}^+} \E\bigl[\|q(t)\|_{\mathcal{I}^+}\bigr] \diff t \\
        &\leq \sqrt{n}\,\widetilde r_{\max}\, \mathcal{B}_q(s) + T \|\theta\|_\infty \sqrt{n}\,\widetilde r_{\max}\, \mathcal{B}_q(s),
    \end{align*}
    where the second inequality follows from the dual bound $\|U^*\|_{\mathcal{I}^+} \leq \sqrt{n}\,\widetilde r_{\max}$ of Lemma~\ref{lmm:dual_bound} in Appendix~\ref{app:dual_bound}, together with the uniform queue bound $\E[\|q(t)\|_{\mathcal{I}^+}] \leq \mathcal{B}_q(s)$ for all $t \geq 0$, established in Proposition~\ref{prop:queue_length_bound}, with
    \begin{align*}
        \mathcal{B}_q(s) = \frac{\sqrt{n}}{\varepsilon} + \left( 1 + \frac{2\sqrt{n}(2e - 3)}{\varepsilon} \right) \ln \left( 2 + \frac{4\sqrt{n}}{\varepsilon} + \frac{4n\|\theta\|_{\mathcal{I}^+}}{s\varepsilon^2} \right).
    \end{align*}
    Dividing both sides by $T$ yields the stated average regret bound:
    \begin{align*}
        F_s^* - \frac{1}{T}\E[F_s^\pi(T)] \leq \left( \frac{1}{T} + \|\theta\|_\infty \right) \sqrt{n}\,\widetilde r_{\max}\, \mathcal{B}_q(s).
    \end{align*}
\end{proof}

\section{Numerical Results}
\label{sec:numerical}
In this section, we numerically evaluate our longest-queue policy against two benchmarks from \citet{aveklouris2021matching}: the matching-rate-based policy and the priority ordering policy. Both benchmarks are discrete-review policies that make matches only at review times, with review period \(\ell_s = \ell s^{-2/3}\) for a constant \(\ell > 0\). The matching-rate-based policy tracks the optimal matching rates of the fluid LP as closely as agent availability allows. The priority ordering policy constructs a fixed ordering of match types based on an optimal extreme-point solution of the fluid LP and greedily makes matches in that order at each review time, subject to agent availability.

We consider three instances of the dynamic matching problem with abandonment. For each instance, we vary the market scale over \(s\in\{10,20,\ldots,200\}\). At each scale, we simulate \(100\) independent replications over a fixed time horizon \(T=100\) to estimate each policy’s approximation ratio. We then compare these ratios to assess how policy performance changes as \(s\) increases.

For these comparisons, we use three review periods for the benchmark policies: \(\ell_s=s^{-2/3}\), \(\ell_s=0.1\), and \(\ell_s=0.3\). The first follows the review rule proposed in \citet{aveklouris2021matching}, while the two fixed periods are also used in their numerical experiments. 
Our longest-queue policy acts at arrival epochs and requires no review
interval, so its performance is unchanged across these choices of
review period.

We begin with two bipartite matching instances each with four demand types, \(D_1,\ldots,D_4\), and four supply types, \(S_1,\ldots,S_4\). The first instance follows the numerical setting in Section~6.1 of \citet{aveklouris2021matching}, with exponential patience times of the same mean and an additional abandonment cost. The demand and supply arrival-rate vectors are \((3,2,1,3)\) and \((2,2,2,2)\), respectively. For each demand type \(D_i\), let \(v_i=(r_{D_iS_1},r_{D_iS_2},r_{D_iS_3},r_{D_iS_4})\) denote its matching-reward vector. We set \(v_1=(1,2,3,1)\), \(v_2=(1,1,1,1)\), \(v_3=(2,1,1,2)\), and \(v_4=(3,3,2,1)\). The holding-cost rates are \((c_{D_1},c_{D_2},c_{D_3},c_{D_4})=(1,2,1,2)\) and \((c_{S_1},c_{S_2},c_{S_3},c_{S_4})=(2,1,2,1)\). For every agent type \(k\), the abandonment rate is \(\theta_k=3\), corresponding to a mean patience time of \(1/3\), and the abandonment cost is \(\gamma_k=1/5\). We also include the terminal cost specified in \eqref{eq:obj}.

Under this parameter setting, the fluid LP has multiple optimal dual solutions, so the first instance violates the GPG condition. For the second instance, we increase only the arrival rate of \(D_1\) from \(3\) to \(4\), keeping all other parameters unchanged. With this change, the optimal fluid dual solution is unique and hence 
the GPG condition holds, i.e., the fluid has a positive GPG. 

In both instances shown in Figures~\ref{fig:numerical-nongpg} and
\ref{fig:numerical-gpg}, our longest-queue policy achieves a higher
approximation ratio than the two benchmark policies for every tested
scale \(s\) and all three review periods \(\ell_s\). The approximation
ratio of our policy approaches \(1\) more rapidly, especially in the second instance,
where the GPG is positive (Figure~\ref{fig:numerical-gpg}). Both benchmark policies
perform worse with longer review periods, as seen when the fixed
interval increases from \(0.1\) to \(0.3\). This is expected because
agents may abandon during the longer wait for the next review, leading
to lost matching opportunities. Our longest-queue policy instead makes
matching decisions at arrival epochs, using a rule that restricts
attention to good matches and gives priority to critical agent types.
It therefore avoids this additional costly delay and outperforms both
benchmark policies in these experiments.

\begingroup
\setlength{\intextsep}{6pt}
\begin{figure}[htbp]
    \centering
    \captionsetup{font=small}
    \captionsetup[subfigure]{font=small}
    \begin{subfigure}[t]{0.325\textwidth}
        \centering
        \includegraphics[width=\linewidth]{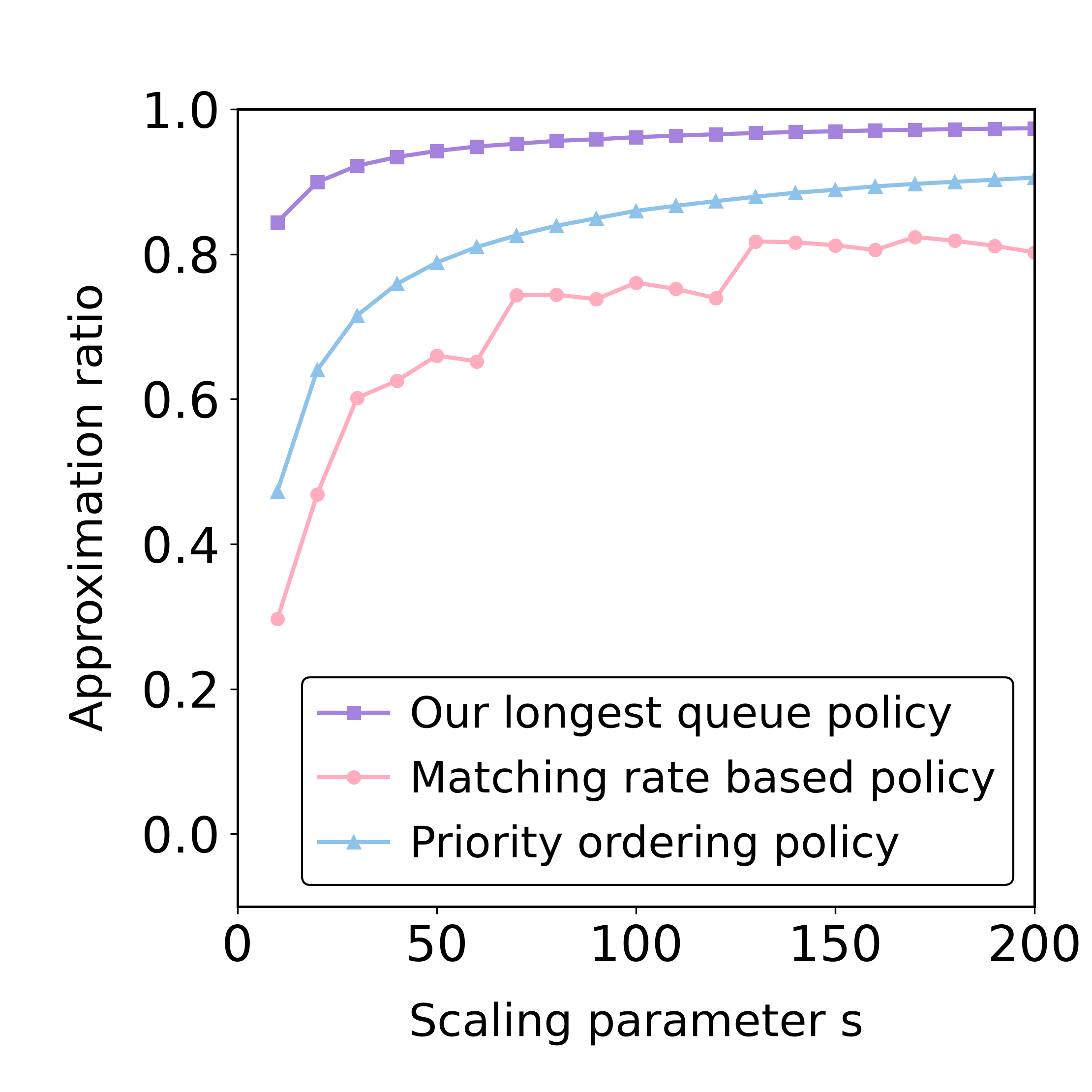}
        \caption{\(\ell_s=s^{-2/3}\)}
        \label{fig:numerical-nongpg-scaled}
    \end{subfigure}\hfill
    \begin{subfigure}[t]{0.325\textwidth}
        \centering
        \includegraphics[width=\linewidth]{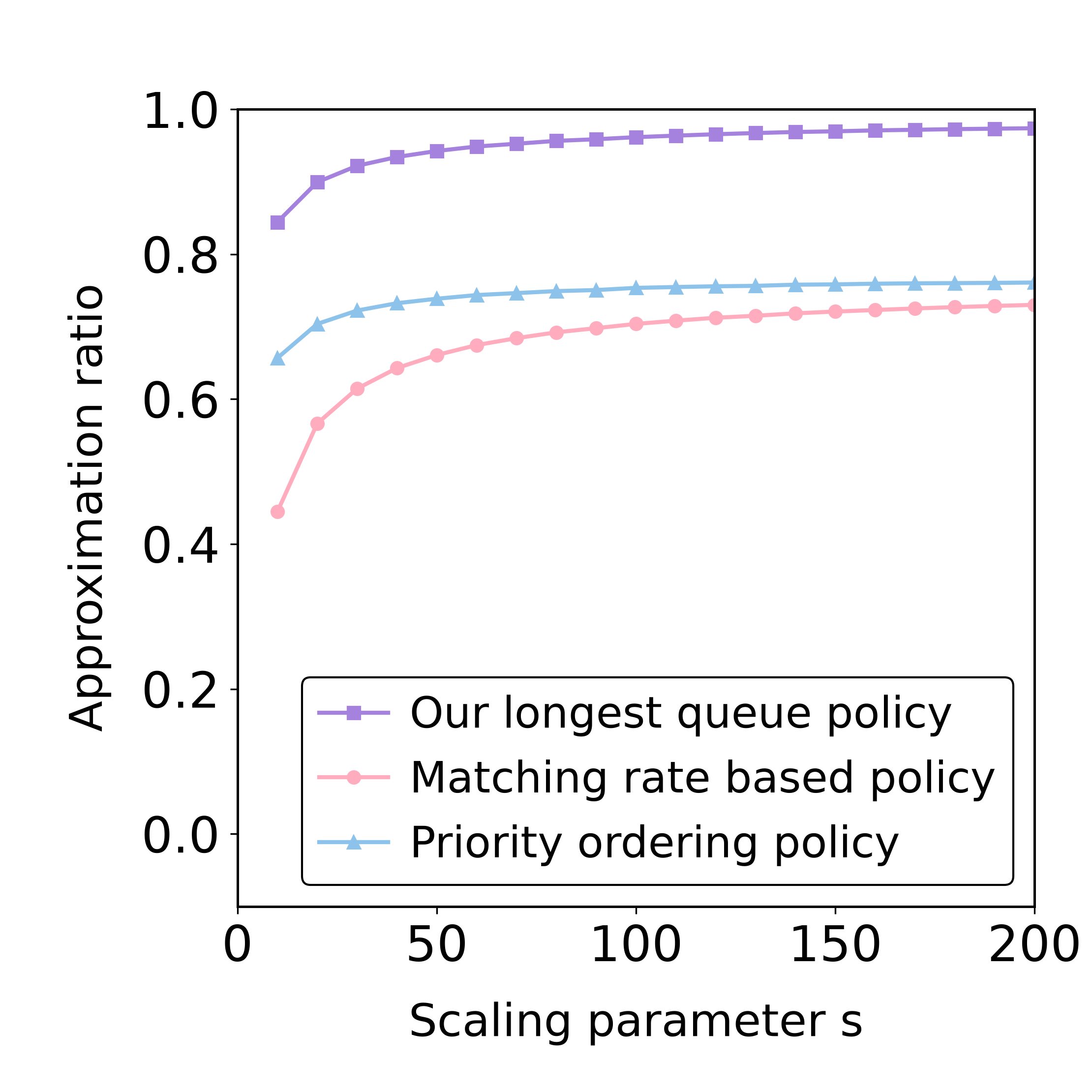}
        \caption{\(\ell_s=0.1\)}
        \label{fig:numerical-nongpg-fixed01}
    \end{subfigure}\hfill
    \begin{subfigure}[t]{0.325\textwidth}
        \centering
        \includegraphics[width=\linewidth]{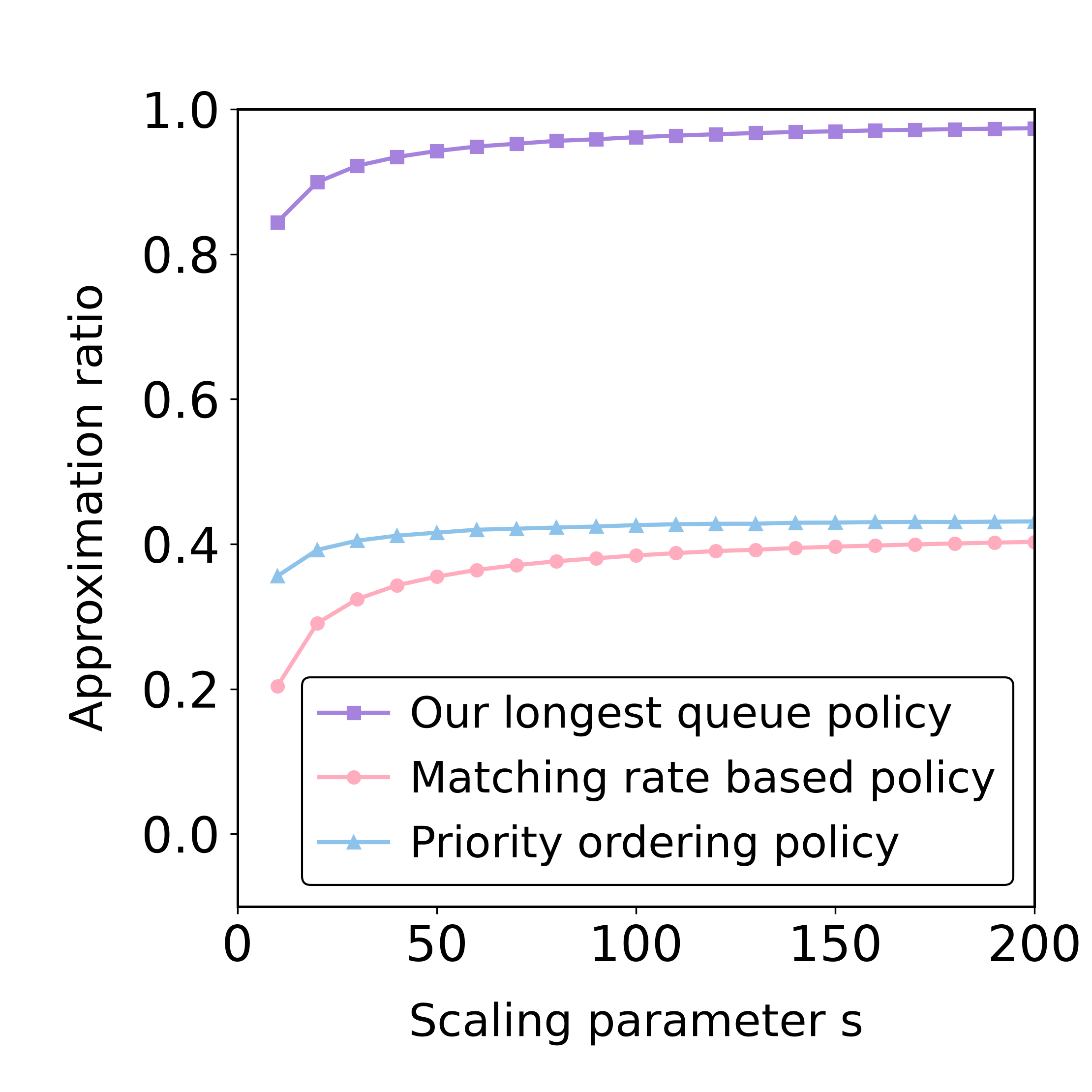}
        \caption{\(\ell_s=0.3\)}
        \label{fig:numerical-nongpg-fixed03}
    \end{subfigure}
    \par\vspace{6pt}
    \caption{The first instance is a bipartite matching network
    with four types on each side, adapted from
    \citet{aveklouris2021matching}. This instance does not satisfy GPG.}
    \label{fig:numerical-nongpg}
\end{figure}

\begin{figure}[htbp]
    \centering
    \captionsetup{font=small}
    \captionsetup[subfigure]{font=small}
    \begin{subfigure}[t]{0.325\textwidth}
        \centering
        \includegraphics[width=\linewidth]{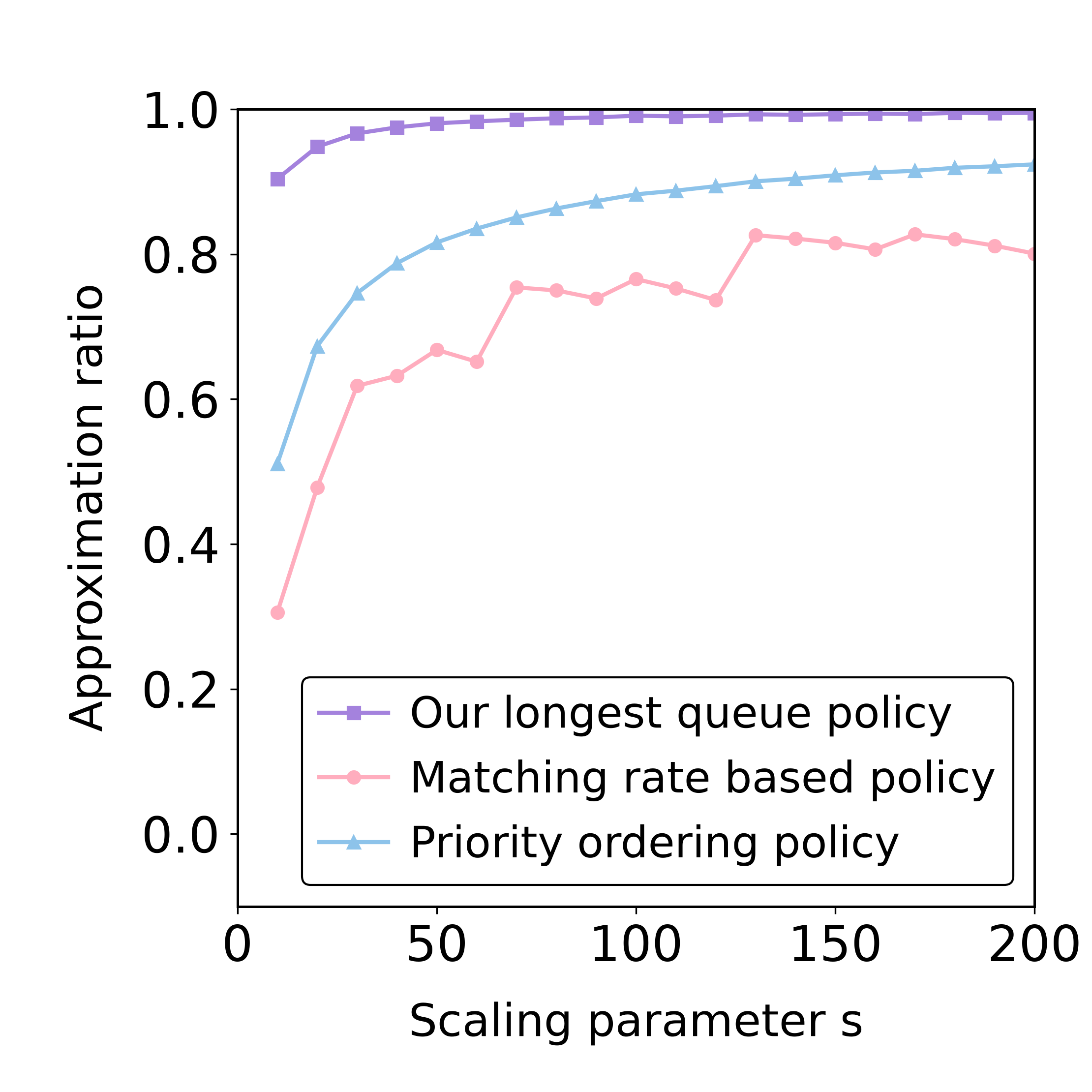}
        \caption{\(\ell_s=s^{-2/3}\)}
        \label{fig:numerical-gpg-scaled}
    \end{subfigure}\hfill
    \begin{subfigure}[t]{0.325\textwidth}
        \centering
        \includegraphics[width=\linewidth]{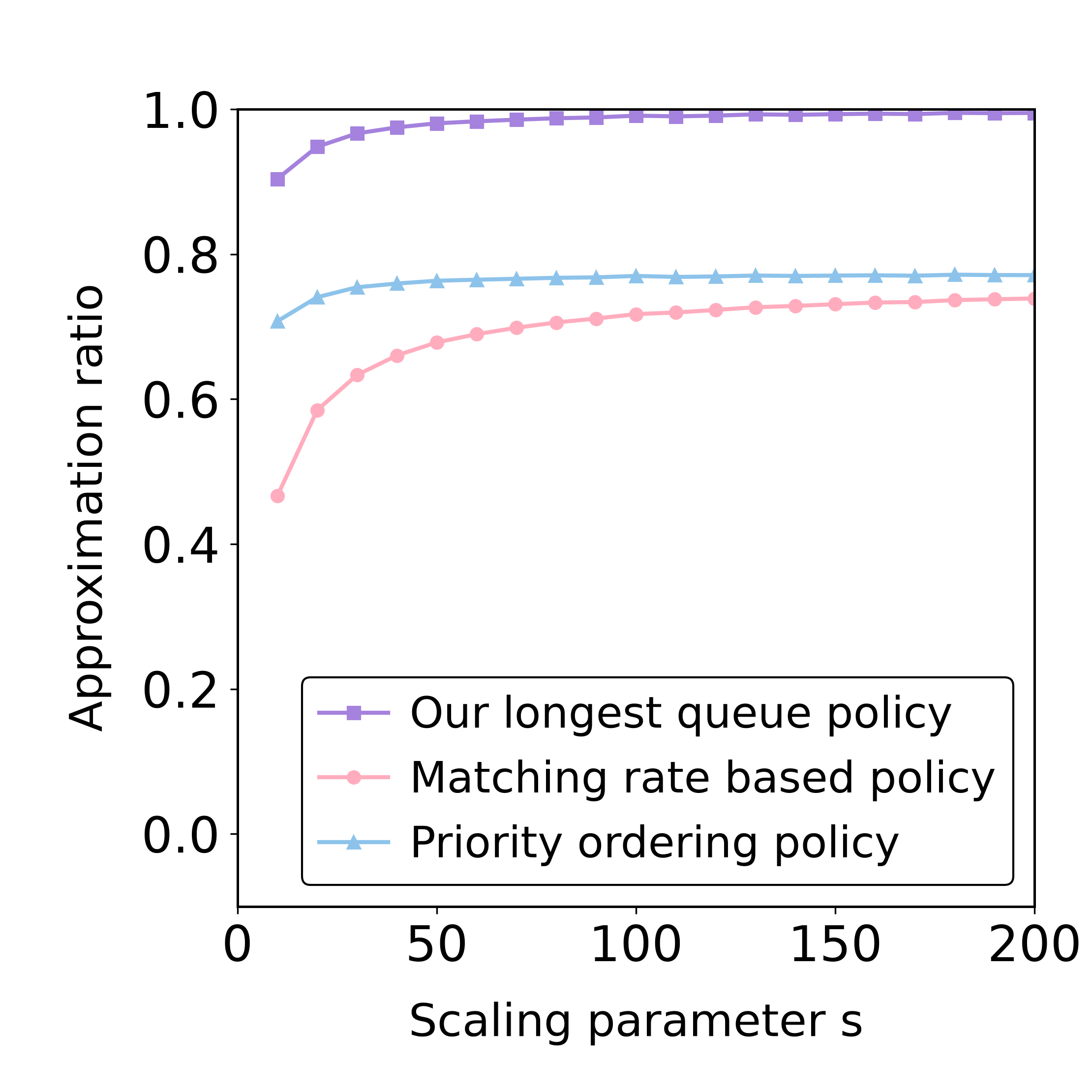}
        \caption{\(\ell_s=0.1\)}
        \label{fig:numerical-gpg-fixed01}
    \end{subfigure}\hfill
    \begin{subfigure}[t]{0.325\textwidth}
        \centering
        \includegraphics[width=\linewidth]{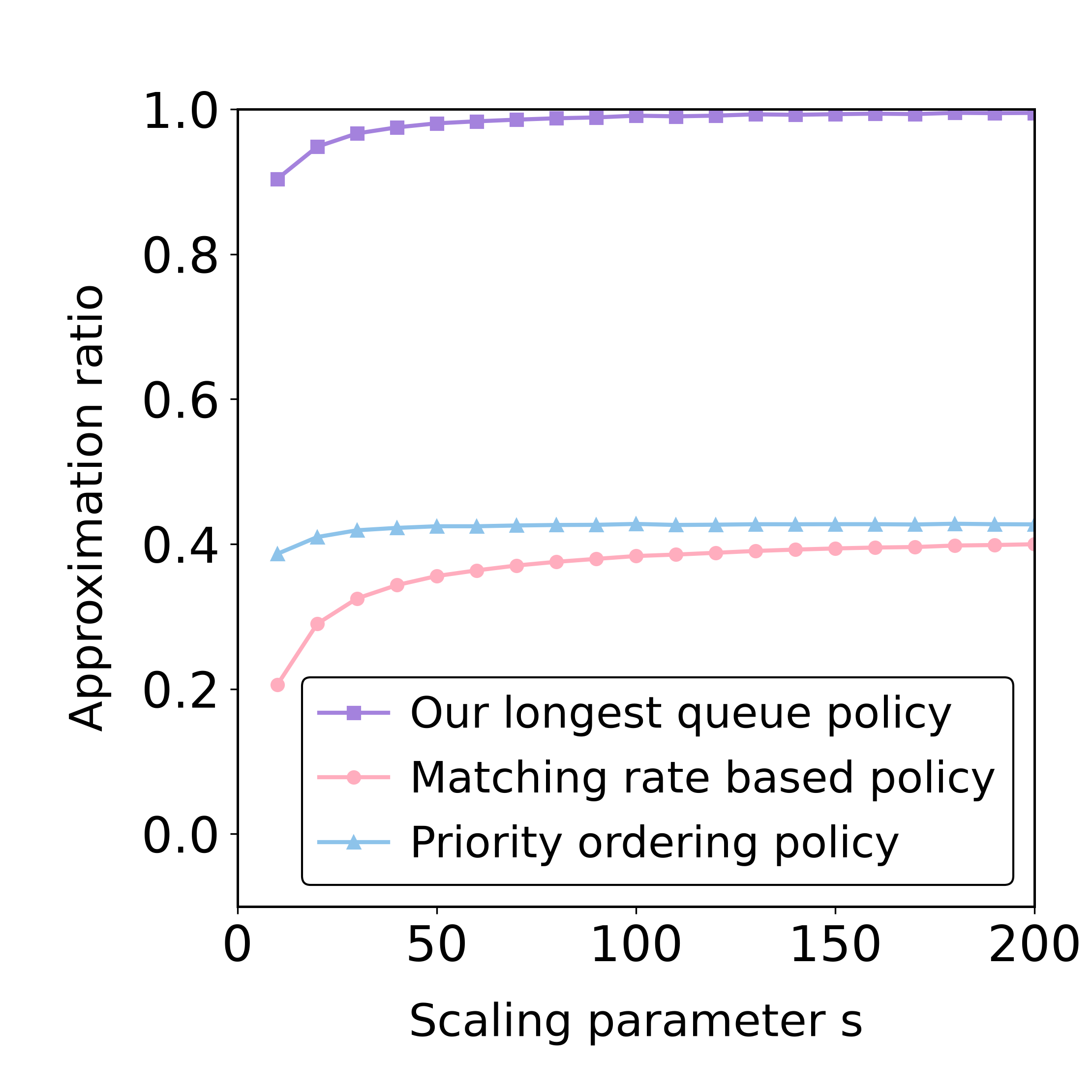}
        \caption{\(\ell_s=0.3\)}
        \label{fig:numerical-gpg-fixed03}
    \end{subfigure}
    \par\vspace{6pt}
    \caption{The second instance is obtained from the first by increasing
    the arrival rate of demand type \(D_1\) from \(3\) to \(4\) and
    satisfies GPG.}
    \label{fig:numerical-gpg}
\end{figure}

\endgroup

\begin{figure}[htbp]
    \centering
    \captionsetup{font=small}

    \begin{subfigure}[t]{0.4\textwidth}
        \centering
        \captionsetup{
            font=small,
            skip=3pt,
            justification=centering,
            singlelinecheck=false,
            margin={0.20\linewidth,0pt}
        }
        \includegraphics[
            width=\linewidth,
            trim={0 24bp 0 0},
            clip
        ]{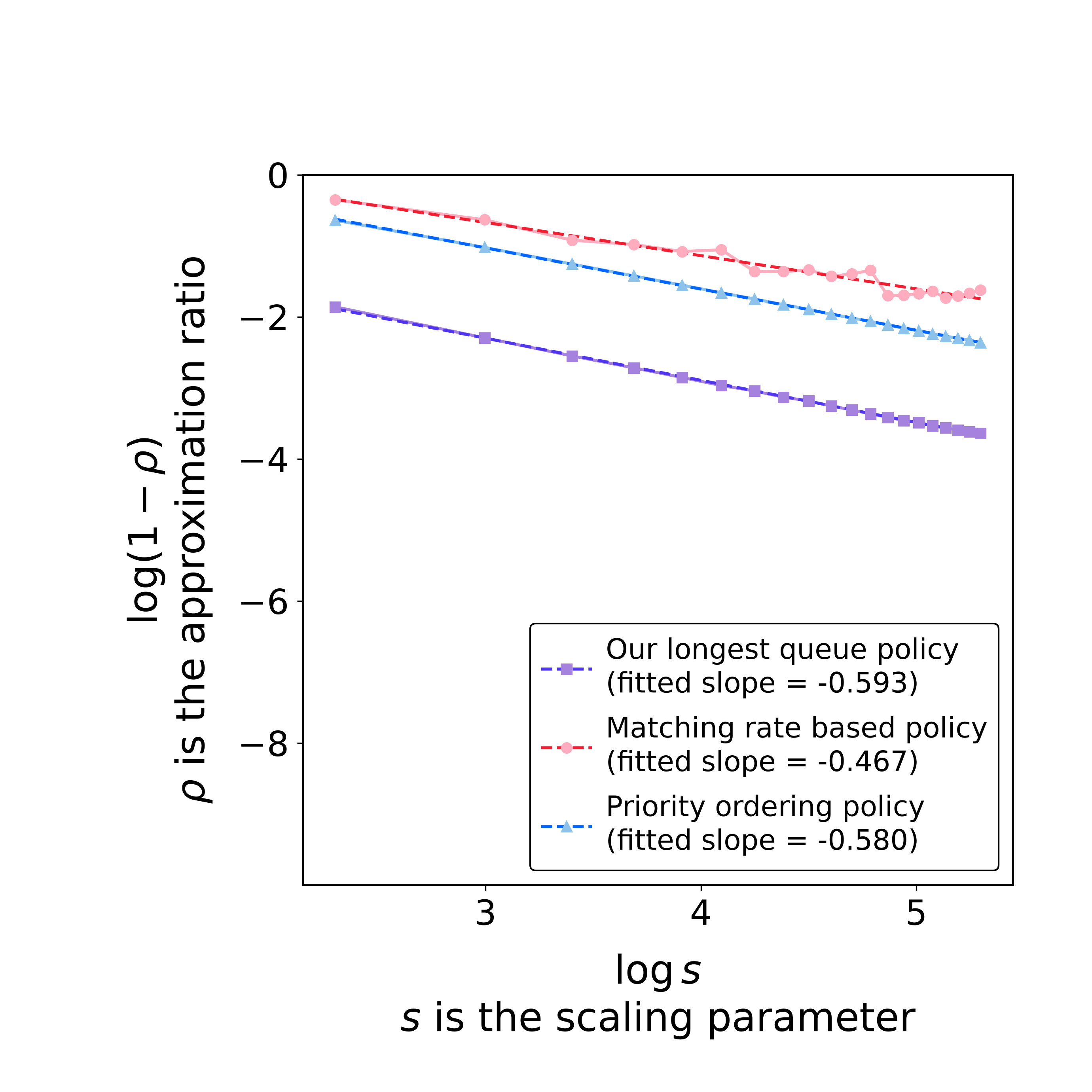}
        \caption{First instance: non-GPG.}
        \label{fig:numerical-loglog-nongpg}
    \end{subfigure}\hspace{0.04\textwidth}
    \begin{subfigure}[t]{0.4\textwidth}
        \centering
        \captionsetup{
            font=small,
            skip=3pt,
            justification=centering,
            singlelinecheck=false,
            margin={0.20\linewidth,0pt}
        }
        \includegraphics[
            width=\linewidth,
            trim={0 24bp 0 0},
            clip
        ]{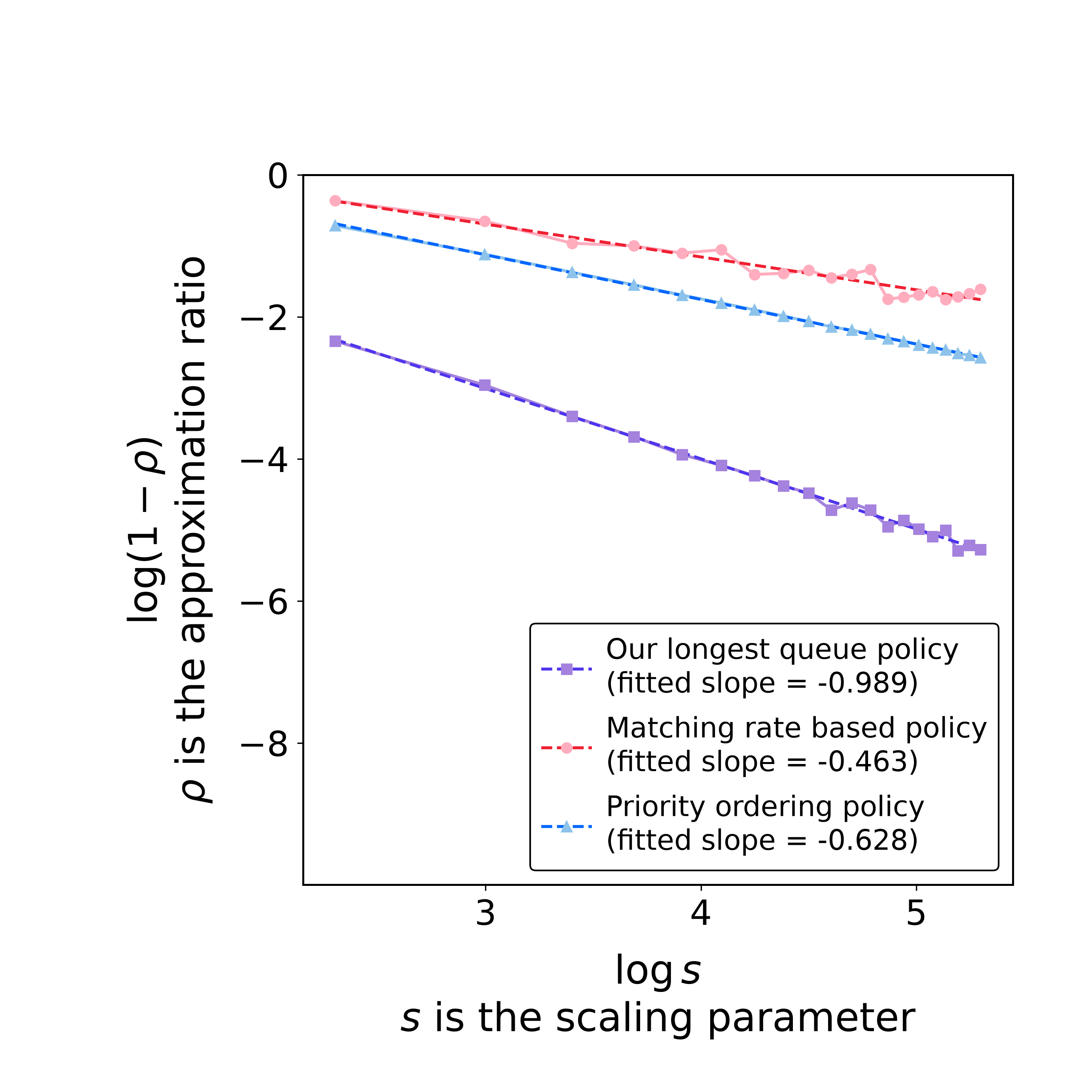}
        \caption{Second instance: GPG.}
        \label{fig:numerical-loglog-gpg}
    \end{subfigure}

    \par\vspace{6pt}
    \caption{\(\log(1-\rho(s))\) versus \(\log s\) for the two bipartite
instances, where \(\rho(s)\) denotes the approximation ratio.
The two benchmark policies use \(\ell_s=s^{-2/3}\).}
\label{fig:numerical-loglog}
\end{figure}

To further compare the policies' convergence rates, we plot
\(\log(1-\rho(s))\) against \(\log s\) in
Figure~\ref{fig:numerical-loglog}, where \(\rho(s)\) denotes the
approximation ratio. For this comparison, we use \(\ell_s=s^{-2/3}\)
for the two benchmark policies, as this choice yields their best
performance at larger market scales. When
\(1-\rho(s)\approx Cs^{-\alpha}\) with \(C,\alpha>0\), taking
logarithms gives
\[
    \log\bigl(1-\rho(s)\bigr)\approx \log C-\alpha\log s.
\]
The slope \(-\alpha\) therefore reflects the convergence rate.

In both instances, our longest-queue policy has the most negative
fitted slope, indicating the fastest convergence rate among the three policies over the
tested scales. In the second instance, where the GPG condition holds, its fitted
slope is very close to \(-1\), consistent with the \(O(1/s)\) bound in
Corollary~\ref{cor:approx_ratio}. In the first instance, where GPG
fails, the fitted slope is close to \(-0.5\), consistent with the
\(O(1/\sqrt{s})\) bound in Proposition~\ref{prop:no_gpg}. 

In the third instance, we consider a non-bipartite two-way matching model
with ten agent types and 24 compatible pairs. The compatibility graph
contains a triangle on types \(1,2,3\) and is therefore non-bipartite.
The base arrival-rate vector is \((7,5,5,3,3,6,2,1,1,3)\).
For every type \(i\), the holding-cost rate is \(c_i=1\), the
abandonment rate is \(\theta_i=3\), and the abandonment cost is
\(\gamma_i=1/5\). Patience times are exponentially distributed with
mean \(1/3\). The rewards for the compatible matching pairs are
\[
    \begin{array}{c|l}
        r_{ij} & \text{Compatible pairs }\{i,j\}\\ \hline
        1 & \{1,7\},\ \{4,5\},\ \{5,6\} \\
        2 & \{1,8\},\ \{6,10\} \\
        3 & \{1,2\},\ \{2,6\},\ \{4,7\},\ \{4,8\},\ \{7,9\},\ \{8,10\},\ \{9,10\} \\
        4 & \{1,3\},\ \{2,4\},\ \{2,7\},\ \{5,9\},\ \{6,9\},\ \{7,8\} \\
        5 & \{2,3\},\ \{2,5\},\ \{3,5\},\ \{3,6\},\ \{3,9\},\ \{5,8\}
    \end{array}
\]
Pairs not listed above are incompatible. 
We note that this instance satisfies the GPG condition. Moreover, while the matching network is non-bipartite, the positive-flow support of the selected optimal
extreme-point solution to the fluid LP is acyclic, which allows us to construct the priority sets following
\citet{aveklouris2021matching} and apply the priority ordering policy.

\begin{figure}[htbp]
    \centering
    \captionsetup{font=small}
    \captionsetup[subfigure]{font=small}
    \begin{subfigure}[t]{0.325\textwidth}
        \centering
        \includegraphics[width=\linewidth]{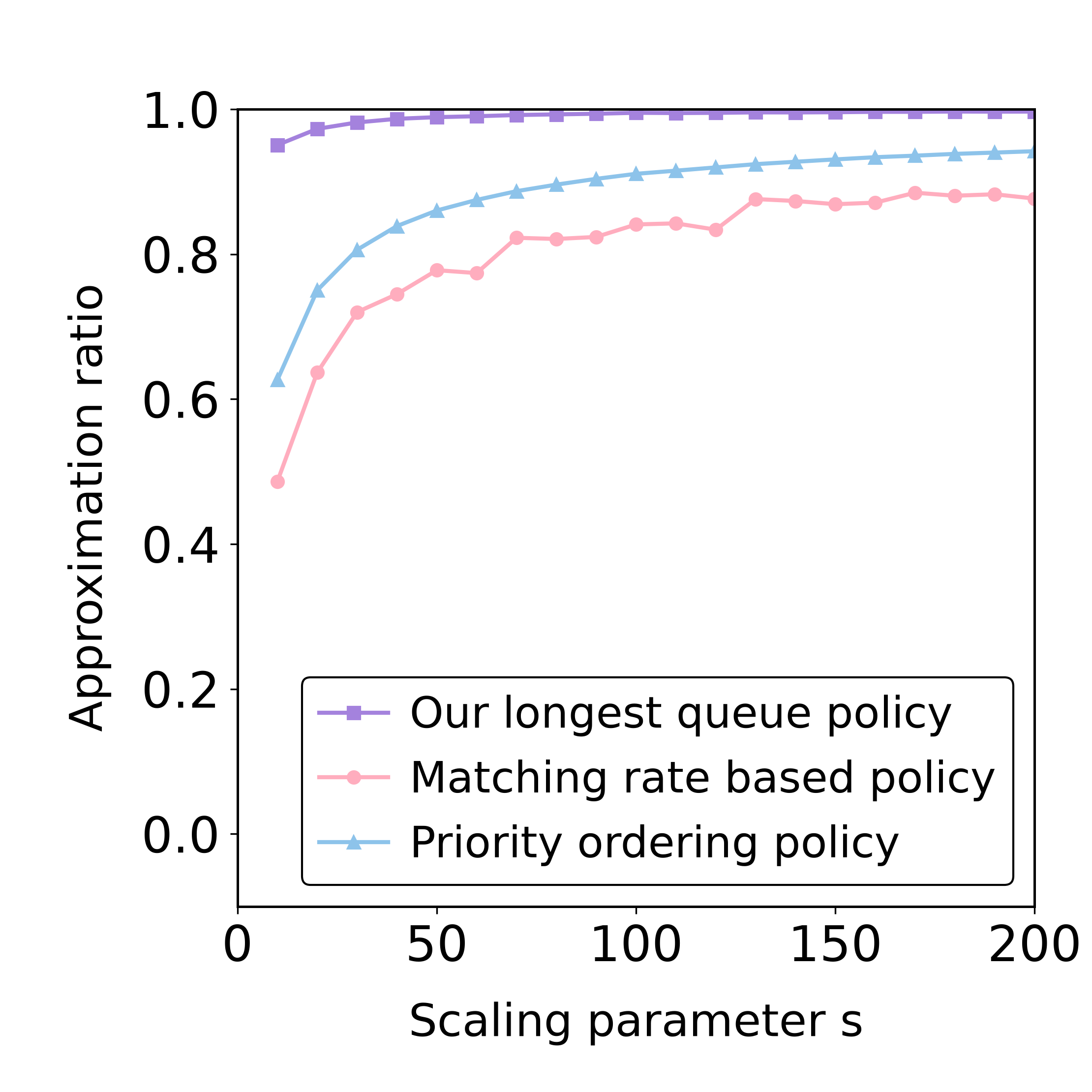}
        \caption{\(\ell_s=s^{-2/3}\)}
        \label{fig:numerical-general-scaled}
    \end{subfigure}\hfill
    \begin{subfigure}[t]{0.325\textwidth}
        \centering
        \includegraphics[width=\linewidth]{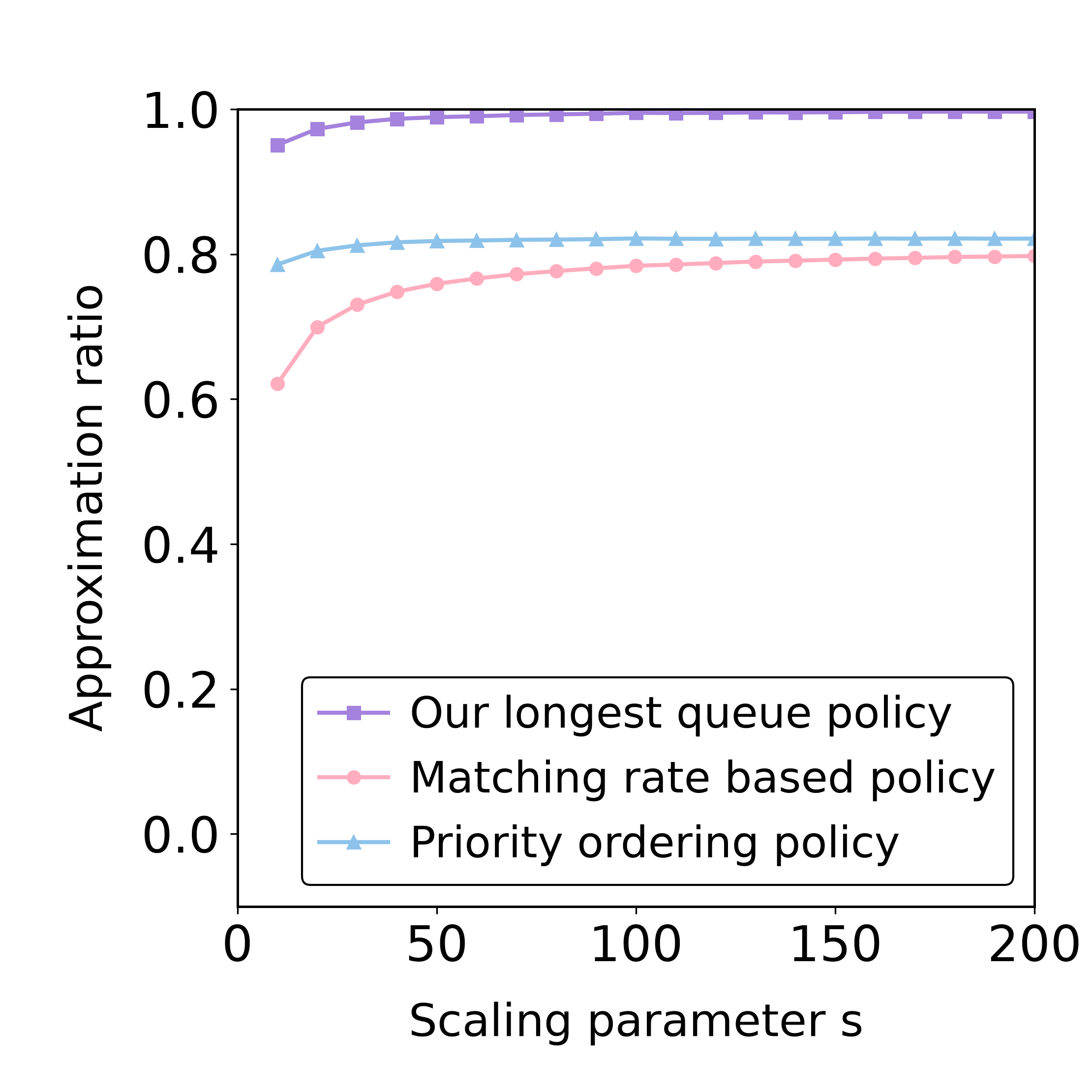}
        \caption{\(\ell_s=0.1\)}
        \label{fig:numerical-general-fixed01}
    \end{subfigure}\hfill
    \begin{subfigure}[t]{0.325\textwidth}
        \centering
        \includegraphics[width=\linewidth]{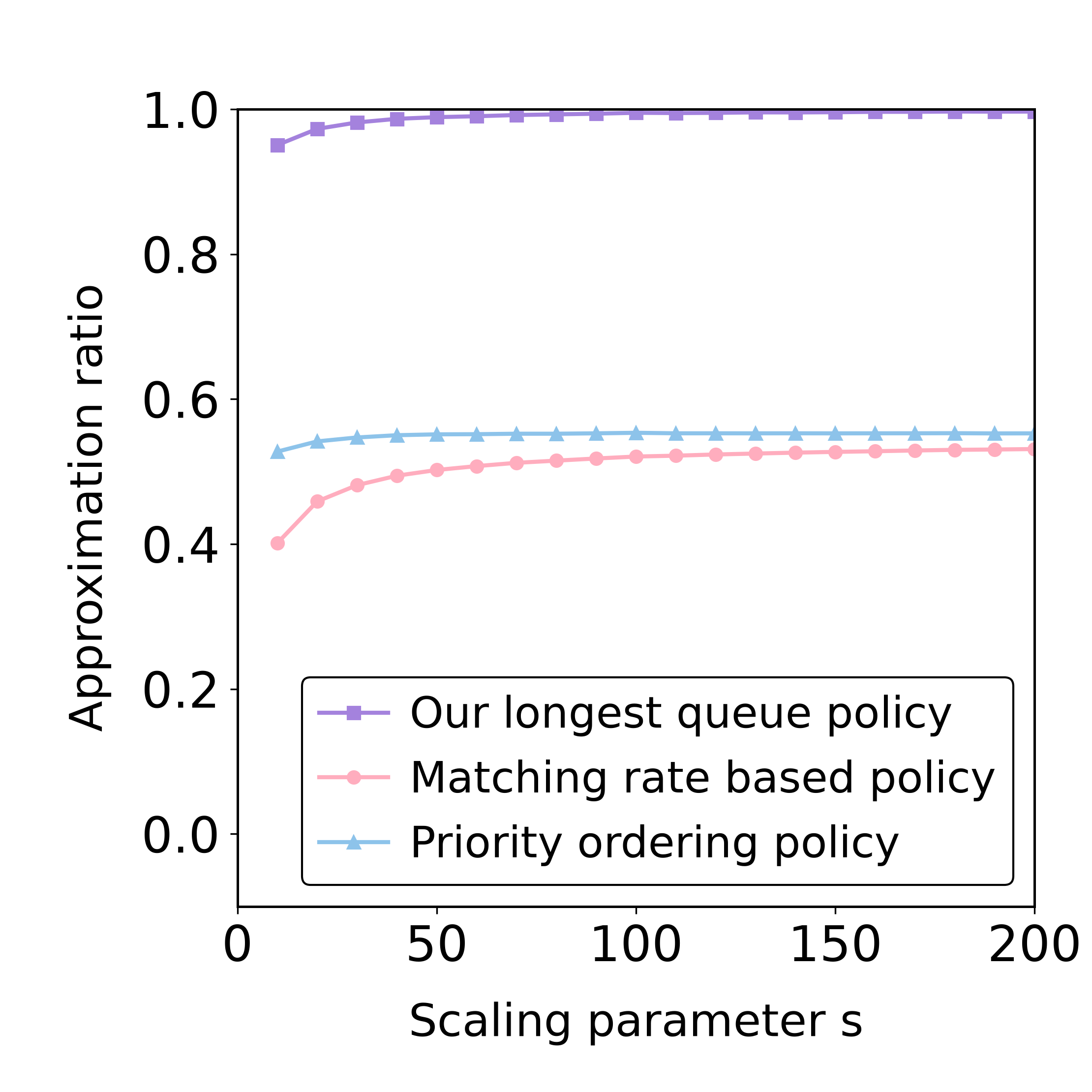}
        \caption{\(\ell_s=0.3\)}
        \label{fig:numerical-general-fixed03}
    \end{subfigure}
    \par\vspace{6pt}
    \caption{The third instance is a non-bipartite two-way matching
    network with ten agent types and \(24\) compatible pairs.}
    \label{fig:numerical-general}
\end{figure}

Figure~\ref{fig:numerical-general} shows that our longest-queue policy
achieves a considerably higher approximation ratio than the two benchmark policies
across all tested scales \(s\) and review periods \(\ell_s\). Its
approximation ratio also approaches \(1\) more rapidly than those of
the two benchmark policies. Thus, the
performance advantage observed in the bipartite instances also appears
in this non-bipartite two-way matching network.

\section{Conclusion}
\label{sec:conclusion}

In this work, we propose a longest-queue policy for dynamic two-way matching
with abandonment. 
Our analysis provides explicit approximation-ratio guarantees for
every \(s\geq1\) and \(T\geq1\). Under exponential patience times,
we prove that our policy achieves an approximation ratio of at least
\(1-C_1/s\) under GPG and at least \(1-C_2/\sqrt{s}\) without GPG,
where \(C_1\) and \(C_2\) are independent of \(s\) and \(T\). The numerical experiments demonstrate that our simple longest-queue policy
outperforms the existing benchmark policies.
The fitted slopes of the log--log plots show that the numerical performance of our policy is consistent with our theoretical bounds.

\section{Methods}
Google Gemini 3 Pro was used through iterative interactions to help develop an
initial proof argument for Lemma~\ref{lmm:general_drift}. The authors
subsequently checked and formalized the argument, and
wrote the proof presented in this paper. OpenAI GPT and Anthropic
Claude models were used to assist in writing the simulation code for the
numerical studies in Section~\ref{sec:numerical}. The authors reviewed and
tested the resulting code and verified the reported numerical results.

\clearpage
\begin{APPENDICES}
\section{Fluid Upper Bound}
\label{app:fluid-upper-bound}

We provide the proof of Lemma~\ref{lmm:fluid-upper-bound} by constructing a feasible
solution to the fluid relaxation whose objective value equals the
policy's expected net value per unit time.
\subsection{Proof of Lemma~\ref{lmm:fluid-upper-bound}}
\begin{proof}[Proof of Lemma~\ref{lmm:fluid-upper-bound}]
Taking expectations in the flow balance
equation~\eqref{eq:flow-balance} and dividing by $T$ yields
\begin{equation}
    \frac{\sum_{m=1}^d M_{im}\E[X_m^\pi(T)]}{T}
    + \frac{\E[q_i^\pi(T)]}{T}
    + \frac{\E[B_i^\pi(T)]}{T}
    = \frac{\E[A_i(T)]}{T}.
    \label{eq:average_flow}
\end{equation}
Since arrivals of type $i$ follow a Poisson process
with rate $\lambda_i$,
\begin{equation}
    \frac{\E[A_i(T)]}{T} = \lambda_i.
    \label{eq:arrival_ev}
\end{equation}

For each $i \in \mathcal{N}$, each waiting type-$i$ agent
abandons according to an independent exponential clock
of rate $\theta_i$. The type-$i$ abandonment intensity
at time $t$ is therefore the number of waiting type-$i$
agents times the individual abandonment rate, namely
$\theta_i q_i^\pi(t)$. Let $Y_i$ denote a unit-rate
Poisson process independent of the arrival processes
and of the other primitive processes. Then the cumulative
number of type-$i$ abandonments admits the
random-time-change representation
\begin{equation*}
    B_i^\pi(t)
    =
    Y_i\!\left(\theta_i\int_0^t q_i^\pi(s)\,ds\right),
    \qquad t \geq 0.
\end{equation*}
The same construction for type-dependent exponential
abandonments is used in
\cite[Section~2.1, equation~(4)]{atar2004scheduling}.

By \cite[Section~3.4, Lemma~3.2]{pang2007martingale},
the compensated time-changed process
\begin{equation}
    M_i^B(t)
    \coloneqq
    B_i^\pi(t)-\theta_i\int_0^t q_i^\pi(s)\,ds
    \label{eq:abandonment-compensator}
\end{equation}
is a martingale. Since $M_i^B(0)=0$, we have
$\mathbb{E}[M_i^B(T)]=0$, and taking expectations
in~\eqref{eq:abandonment-compensator} gives
\begin{equation}
    \mathbb{E}[B_i^\pi(T)]
    =
    \theta_i\mathbb{E}\!\left[
        \int_0^T q_i^\pi(t)\,dt
    \right]
    =
    \theta_i\int_0^T\mathbb{E}[q_i^\pi(t)]\,dt,
    \label{eq:expected-abandonments}
\end{equation}
where the last equality follows from Fubini--Tonelli's
theorem \cite[Section~2.5, Theorem~2.37]{folland1999real},
using $q_i^\pi \geq 0$.

Define the time-averaged queue length and matching rates as
\begin{equation}
    \bar q_i
    \triangleq
    \frac{1}{T}\int_0^T\E[q_i^\pi(t)]\,dt
    + \frac{\E[q_i^\pi(T)]}{T\theta_i},
    \qquad
    \bar x_m
    \triangleq
    \frac{\E[X_m^\pi(T)]}{T}.
    \label{eq:avg_queue}
\end{equation}
Combining~\eqref{eq:expected-abandonments}
with~\eqref{eq:avg_queue}, we have
\begin{equation}
    \theta_i\bar q_i
    =
    \frac{\E[B_i^\pi(T)]}{T}
    + \frac{\E[q_i^\pi(T)]}{T}.
    \label{eq:queue_relation}
\end{equation}
Substituting~\eqref{eq:arrival_ev},~\eqref{eq:avg_queue},
and~\eqref{eq:queue_relation}
into~\eqref{eq:average_flow} gives
\begin{equation*}
    [M\bar x]_i+\theta_i\bar q_i=\lambda_i,
\end{equation*}
which is the equality constraint in~\eqref{eq:fluidwithq}.
Since $\bar x_m,\bar q_i\geq0$ by construction,
the pair $(\bar x,\bar q)$ is feasible
for~\eqref{eq:fluidwithq}.

It remains to show that $\mathbb{E}[F^\pi(T)]/T$
equals the objective of~\eqref{eq:fluidwithq}
evaluated at $(\bar x,\bar q)$.
Taking expectations of $F^\pi(T)$ as defined
in~\eqref{eq:obj} and dividing by $T$, we have
\begin{align*}
    \frac{\E[F^\pi(T)]}{T}
    &=
    \sum_{m=1}^d r_m\frac{\E[X_m^\pi(T)]}{T}
    - \sum_{i=1}^n c_i\cdot
      \frac{1}{T}\int_0^T\E[q_i^\pi(t)]\,dt\\
    &\quad
    - \sum_{i=1}^n\gamma_i\frac{\E[B_i^\pi(T)]}{T}
    - \sum_{i=1}^n
      \left(\frac{c_i}{\theta_i}+\gamma_i\right)
      \frac{\E[q_i^\pi(T)]}{T}\\[6pt]
    &=
    \sum_{m=1}^d r_m\bar x_m
    - \sum_{i=1}^n c_i
      \underbrace{
          \left(
              \frac{1}{T}\int_0^T\E[q_i^\pi(t)]\,dt
              + \frac{\E[q_i^\pi(T)]}{T\theta_i}
          \right)
      }_{=\bar q_i}\\
    &\quad
    - \sum_{i=1}^n\gamma_i
      \underbrace{
          \left(
              \frac{\E[B_i^\pi(T)]}{T}
              + \frac{\E[q_i^\pi(T)]}{T}
          \right)
      }_{=\theta_i\bar q_i}\\[6pt]
    &=
    \sum_{m=1}^d r_m\bar x_m
    - \sum_{i=1}^n(c_i+\theta_i\gamma_i)\bar q_i\\
    &=
    \sum_{m=1}^d r_m\bar x_m
    - \sum_{i=1}^n\bar c_i\bar q_i,
\end{align*}
where we used~\eqref{eq:avg_queue}
and~\eqref{eq:queue_relation}.
This is exactly the objective of~\eqref{eq:fluidwithq}
evaluated at $(\bar x,\bar q)$.
Since $F^*(\lambda)$ is the optimal value
of~\eqref{eq:fluidwithq}, we conclude
\begin{equation*}
    \frac{\mathbb E[F^\pi(T)]}{T}
    \leq F^*(\lambda).
\end{equation*}
\end{proof}

\section{Supplementary Proofs for Section~\ref{sec:analysis-lq}}
\label{app:analysis-lq}

This section provides the proofs deferred from
Section~\ref{sec:analysis-lq}.

\subsection{Proofs for Section~\ref{sec:step1}}
\label{app:step1}

\subsubsection{Proof of Lemma~\ref{lmm:bound_reward}}\label{app:lmm:bound_reward}
\begin{proof}[Proof of Lemma~\ref{lmm:bound_reward}]
Taking expectations of $F^\pi(T)$ as defined
in~\eqref{eq:obj} yields
\begin{align}
    \E[F^\pi(T)]
    &=
    \E\left[\sum_{m=1}^{d} r_m X_m^\pi(T)\right]
    - \sum_{i=1}^{n} c_i
      \E\left[\int_0^T q_i^\pi(t)\,dt\right]
    \notag\\
    &\quad
    - \sum_{i=1}^{n} \gamma_i \E[B_i^\pi(T)]
    - \sum_{i=1}^{n}
      \left(\frac{c_i}{\theta_i}+\gamma_i\right)
      \E[q_i^\pi(T)].
    \label{eq:expected_obj}
\end{align}
As established in~\eqref{eq:expected-abandonments},
the expected cumulative number of type-$i$
abandonments satisfies
\begin{equation*}
    \E[B_i^\pi(T)]
    =
    \theta_i\E\left[\int_0^T q_i^\pi(t)\,dt\right].
\end{equation*}
Substituting this identity into the expectation of the
flow-balance equation~\eqref{eq:flow-balance}
and rearranging yields
\begin{align}
    \E\left[\int_0^T q_i^\pi(t)\,dt\right]
    =
    \frac{1}{\theta_i}
    \left(
        \lambda_i T
        - \E[q_i^\pi(T)]
        - \E[(MX^\pi(T))_i]
    \right),
    \label{eq:cum_queue}
\end{align}
where we used $\E[A_i(T)]=\lambda_i T$ since arrivals
follow a Poisson process with rate $\lambda_i$.

Substituting~\eqref{eq:expected-abandonments}
and~\eqref{eq:cum_queue} into~\eqref{eq:expected_obj} gives
\begin{align*}
    \E[F^\pi(T)]
    &=
    \E\left[\sum_{m=1}^{d} r_m X_m^\pi(T)\right]
    - \sum_{i=1}^{n}\frac{c_i}{\theta_i}
      \left(
          \lambda_i T
          - \E[q_i^\pi(T)]
          - \E[(MX^\pi(T))_i]
      \right)\\
    &\quad
    - \sum_{i=1}^{n}\gamma_i
      \left(
          \lambda_i T
          - \E[q_i^\pi(T)]
          - \E[(MX^\pi(T))_i]
      \right)
    - \sum_{i=1}^{n}
      \left(\frac{c_i}{\theta_i}+\gamma_i\right)
      \E[q_i^\pi(T)]\\[6pt]
    &=
    \E\left[\sum_{m=1}^{d} r_m X_m^\pi(T)\right]
    - \sum_{i=1}^{n}\frac{\bar c_i}{\theta_i}
      \left(
          \lambda_i T
          - \E[q_i^\pi(T)]
          - \E[(MX^\pi(T))_i]
      \right)
    - \sum_{i=1}^{n}\frac{\bar c_i}{\theta_i}
      \E[q_i^\pi(T)],
\end{align*}
where the second equality follows from the definition
of $\bar c_i$ in~\eqref{eq:adjusted-holding-cost}.
The $\E[q_i^\pi(T)]$ terms cancel, yielding
\begin{align*}
    \E[F^\pi(T)]
    &=
    \E\left[\sum_{m=1}^{d} r_m X_m^\pi(T)\right]
    + \sum_{i=1}^{n}\frac{\bar c_i}{\theta_i}
      \E[(MX^\pi(T))_i]
    - \sum_{i=1}^{n}\frac{\bar c_i\lambda_i T}{\theta_i}
    \\[6pt]
    &=
    \E\left[
        \sum_{m=1}^{d}
        \left(
            r_m+\sum_{i=1}^{n}
            \frac{\bar c_i}{\theta_i}M_{im}
        \right)X_m^\pi(T)
    \right]
    - \sum_{i=1}^{n}\frac{\bar c_i\lambda_i T}{\theta_i}
    \\[6pt]
    &=
    \E\left[\sum_{m=1}^{d}\tilde r_m X_m^\pi(T)\right]
    - \sum_{i=1}^{n}\frac{\bar c_i\lambda_i T}{\theta_i},
\end{align*}
where the last equality follows from the definition
of $\tilde r_m$ in~\eqref{eq:adjusted-match-reward}.
Equation~\eqref{eq:fluid-value-at-x-star} gives the
optimal fluid objective value $F^*$; multiplying it
by $T$ gives
\begin{equation*}
    T F^*
    =
    T\langle\tilde r,x^*\rangle
    - \sum_{i=1}^{n}\frac{\bar c_i\lambda_i T}{\theta_i}.
\end{equation*}
Subtracting the expression for $\E[F^\pi(T)]$ from $TF^*$
cancels the common constant, yielding
\begin{align*}
    T F^*-\E[F^\pi(T)]
    =
    T\langle\tilde r,x^*\rangle
    - \E\left[\sum_{m=1}^{d}\tilde r_m X_m^\pi(T)\right].
\end{align*}
\end{proof}

\subsection{Proofs for Section~\ref{sec:step2}}
\label{app:step2}

We provide the proofs of the auxiliary lemmas used to
establish Proposition~\ref{prop:queue_square_negdrift}.

\subsubsection{Proof of Lemma~\ref{lmm:generator_bound}}
\label{app:lmm:generator_bound}
Using Definition~\ref{def:generator}, we first derive an
expression for $\mathcal L f(q(t))$ under the
longest-queue policy.

\begin{lemma}
\label{lmm:generator}
Under Algorithm~\ref{alg:shadow-matching-neighbor-set}, the queue length process $\{q(t)\}_{t \ge 0}$ evolves as a time-homogeneous continuous-time Markov chain on $\mathbb{Z}_+^n$. Given the current state $q(t) = q$, the process is driven by two types of events:
\begin{enumerate}
    \item \textbf{Arrivals:} Type-$i$ agents arrive as a Poisson process with rate $\lambda_i$, causing $q \to q + e_i - Mx^{i}(t)$, as in~\eqref{eq:fg-queue-update-arrival}.
    \item \textbf{Abandonments:} Each agent in queue $i$ abandons independently at rate $\theta_i$, causing $q \to q - e_i$, as in~\eqref{eq:fg-queue-update-abandonment}.
\end{enumerate}
The infinitesimal generator acting on a bounded measurable function $f$ is
\begin{align}\label{eq:generator}
    \mathcal{L} f(q) = \sum_{i=1}^n \lambda_i \bigl[ f(q + e_i - Mx^{i}(t)) - f(q) \bigr] + \sum_{i=1}^n \theta_i q_i \bigl[ f(q - e_i) - f(q) \bigr].
\end{align}
\end{lemma}

\begin{proof}[Proof of Lemma~\ref{lmm:generator}]
By Definition~\ref{def:generator}, the generator is $\mathcal{L} f(q) = \lim_{h \downarrow 0} \frac{\mathbb{E}[f(q(t+h)) \mid q(t) = q] - f(q)}{h}$. We compute this limit by analyzing the event probabilities over $[t, t+h)$. By the standard infinitesimal transition probabilities for independent Poisson arrivals and exponential abandonments, we have
\begin{equation}
\begin{aligned}
    \mathbb{P}(\text{one type-}i \text{ arrival in } [t,t+h) \mid q(t) = q) &= \lambda_i h + o(h), \\
    \mathbb{P}(\text{one type-}i \text{ abandonment in } [t,t+h) \mid q(t) = q) &= \theta_i q_i h + o(h), \\
    \mathbb{P}(\text{two or more events in } [t,t+h) \mid q(t) = q) &= o(h), \\
    \mathbb{P}(\text{no event in } [t,t+h) \mid q(t) = q) &= 1 - \Lambda(q) h + o(h),
\end{aligned}
\label{eq:event-probabilities}
\end{equation}
where $\Lambda(q) = \sum_{i=1}^n \lambda_i + \sum_{i=1}^n \theta_i q_i$ is the
total event rate, and $g(h) = o(h)$ means
$\lim_{h \downarrow 0} g(h)/h = 0$.  

Conditioning on which event occurs and using the event probabilities~\eqref{eq:event-probabilities} together with the transitions~\eqref{eq:fg-queue-update-arrival} and~\eqref{eq:fg-queue-update-abandonment}, we have
\begin{align*}
    \mathbb{E}[f(q(t+h)) \mid q(t) = q]
    &= f(q)(1 - \Lambda(q) h) + \sum_{i=1}^n \lambda_i h \cdot f(q + e_i - Mx^{i}(t)) \\
    &\quad + \sum_{i=1}^n \theta_i q_i h \cdot f(q - e_i) + o(h).
\end{align*}
Subtracting $f(q)$ and dividing by $h$ yields
\begin{align*}
    \frac{\mathbb{E}[f(q(t+h)) \mid q(t) = q] - f(q)}{h}
    = \sum_{i=1}^n \lambda_i \bigl[ f(q + e_i - Mx^{i}(t)) - f(q) \bigr] + \sum_{i=1}^n \theta_i q_i \bigl[ f(q - e_i) - f(q) \bigr] + \frac{o(h)}{h}.
\end{align*}
Taking $h \downarrow 0$ yields~\eqref{eq:generator}.
\end{proof}

Having established the generator in Lemma~\ref{lmm:generator},
we now prove Lemma~\ref{lmm:generator_bound} by applying it to
the squared norm of the critical queues.

\begin{proof}[Proof of Lemma~\ref{lmm:generator_bound}]
Applying the generator~\eqref{eq:generator} from Lemma~\ref{lmm:generator} to $f(q) = \|q\|^2_{\mathcal{I}^+} = \sum_{i \in \mathcal{I}^+} q_i^2$ yields
\begin{align*}
    \mathcal{L}\|q(t)\|^2_{\mathcal{I}^+}
    &= \sum_{i \in \mathcal{I}^+} \lambda_i^s \left[ \left( q_i(t) + 1 - [Mx^{i}(t)]_i \right)^2 + \sum_{j \neq i, j \in \mathcal{I}^+} \left( q_j(t) - [Mx^{i}(t)]_j \right)^2 - \sum_{j \in \mathcal{I}^+} q_j(t)^2 \right] \\
    &\quad + \sum_{i \in [n] \setminus \mathcal{I}^+} \lambda_i^s \left[ \sum_{j \in \mathcal{I}^+} \left( q_j(t) - [Mx^{i}(t)]_j \right)^2 - \sum_{j \in \mathcal{I}^+} q_j(t)^2 \right] \\
    &\quad + \sum_{i \in \mathcal{I}^+} \theta_i q_i(t) \left[ \left( q_i(t) - 1 \right)^2 - q_i(t)^2 \right],
\end{align*}
where the abandonment sum is restricted to $i \in \mathcal{I}^+$ since abandonments from $i \notin \mathcal{I}^+$ do not affect $\|q(t)\|^2_{\mathcal{I}^+}$.
Expanding the squared terms gives
\begin{align}
    \mathcal{L}\|q(t)\|^2_{\mathcal{I}^+}
    &= \sum_{i \in \mathcal{I}^+} \lambda_i^s
    \left[
        2q_i(t)
        - 2\langle q(t), Mx^{i}(t) \rangle_{\mathcal{I}^+}
        + 1 - 2[Mx^{i}(t)]_i
        + \|Mx^{i}(t)\|^2_{\mathcal{I}^+}
    \right] \notag\\
    &\quad
    + \sum_{i \in [n] \setminus \mathcal{I}^+} \lambda_i^s
    \left[
        -2\langle q(t), Mx^{i}(t) \rangle_{\mathcal{I}^+}
        + \|Mx^{i}(t)\|^2_{\mathcal{I}^+}
    \right] \notag\\
    &\quad
    + \sum_{i \in \mathcal{I}^+}
    \theta_i q_i(t)\left[-2q_i(t)+1\right].
    \label{eq:drift-expanded-abandonment}
\end{align}
Since $q_i(t)$ is a nonnegative integer, $q_i(t)\left[ -2q_i(t) + 1 \right] \leq 0$, so the abandonment term is non-positive and we drop it. Combining the remaining terms and substituting the definition~\eqref{eq:aggregate-matching-rate} of $\bar{x}(t)$, we obtain
\begin{align}
    \mathcal{L}\|q(t)\|^2_{\mathcal{I}^+}
    \leq 2\langle q(t), \lambda^s - M\bar{x}(t) \rangle_{\mathcal{I}^+} + \sum_{i \in \mathcal{I}^+} \lambda_i^s \left( 1 - 2[Mx^{i}(t)]_i \right) + \sum_{i =1}^n \lambda_i^s \|Mx^{i}(t)\|^2_{\mathcal{I}^+}.
    \label{eq:drift-combined}
\end{align}

It remains to bound the last two sums. We combine them into a single sum,
\begin{align*}
    \sum_{i \in \mathcal{I}^+} \lambda_i^s \left( 1 - 2[Mx^{i}(t)]_i \right) + \sum_{i=1}^n \lambda_i^s \|Mx^{i}(t)\|^2_{\mathcal{I}^+}
    = \sum_{i=1}^n \lambda_i^s \kappa_i(t),
\end{align*}
where
\begin{align*}
    \kappa_i(t) \coloneqq \mathbf{1}\{i \in \mathcal{I}^+\}\left( 1 - 2[Mx^{i}(t)]_i \right) + \|Mx^{i}(t)\|^2_{\mathcal{I}^+}.
\end{align*}
 We show that $\kappa_i(t) \leq 1$ for every $i$. Since $Mx^{i}(t)$ has entries in $\{0,1\}$, we have
\begin{align*}
    \|Mx^{i}(t)\|^2_{\mathcal{I}^+} = \sum_{k \in \mathcal{I}^+} \left( [Mx^{i}(t)]_k \right)^2 = \sum_{k \in \mathcal{I}^+} [Mx^{i}(t)]_k.
\end{align*}
Upon a type-$i$ arrival, either no match is selected or the arriving agent is matched with exactly one waiting agent.

If no match is selected, then $x^{i}(t) = 0$ and $Mx^{i}(t) = 0$, so $\kappa_i(t) = \mathbf{1}\{i \in \mathcal{I}^+\} \leq 1$.

If a match $m$ with $\mathcal{N}(m) = \{i,j\}$ is selected, then $Mx^{i}(t) = e_i + e_j$, so $[Mx^{i}(t)]_i = 1$ and $\|Mx^{i}(t)\|^2_{\mathcal{I}^+} = \mathbf{1}\{i \in \mathcal{I}^+\} + \mathbf{1}\{j \in \mathcal{I}^+\}$. Hence
\begin{align*}
    \kappa_i(t) = -\mathbf{1}\{i \in \mathcal{I}^+\} + \mathbf{1}\{i \in \mathcal{I}^+\} + \mathbf{1}\{j \in \mathcal{I}^+\} = \mathbf{1}\{j \in \mathcal{I}^+\} \leq 1.
\end{align*}

In both cases $\kappa_i(t) \leq 1$, so
\begin{align}
\sum_{i \in \mathcal{I}^+} \lambda_i^s \left( 1 - 2[Mx^{i}(t)]_i \right) + \sum_{i=1}^n \lambda_i^s \|Mx^{i}(t)\|^2_{\mathcal{I}^+} \leq \sum_{i=1}^n \lambda_i^s = \|\lambda^s\|_1 = s,
    \label{eq:jump-bound}
\end{align}
where the last equality uses $\|\lambda^s\|_1 = s\|\lambda^b\|_1 = s$. Combining~\eqref{eq:drift-combined} and~\eqref{eq:jump-bound}, we obtain
\begin{align*}
    \mathcal{L}\|q(t)\|^2_{\mathcal{I}^+} \leq 2s\langle q(t), \lambda^b - M\bar{x}^b(t) \rangle_{\mathcal{I}^+} + s.
\end{align*}
\end{proof}

\subsubsection{Proof of Lemma~\ref{lmm:optimal_action}}
\label{app:lmm:optimal_action}
Before proving the optimality result in Lemma~\ref{lmm:optimal_action},
we first establish a simple property of our algorithm.
Because the policy is maximally aggressive on good matches, no two waiting agents can ever form a good match. Whenever a new agent arrives, the algorithm immediately matches it if some waiting agent forms a good match with it. Hence if such a pair of waiting agents existed, the later arrival would have been matched on arrival, since the earlier agent was already present and available. 

\begin{lemma}\label{lmm:empty_neighbor}
    Under Algorithm~\ref{alg:shadow-matching-neighbor-set}, for every good match $(i,j) \in \mathcal{M}^g$ and every $t \geq 0$,
    \begin{align*}
        q_i(t)\, q_j(t) = 0.
    \end{align*}
    In particular, if $i \in \mathcal{Q}^+(t)$, then $\mathcal{N}^+_i(t) = \mathcal{N}^0_i(t) = \emptyset$, where $\mathcal N_i^+(t)$ and $\mathcal N_i^0(t)$
are defined in~\eqref{eq:fg-available-critical-neighbors}
and~\eqref{eq:fg-available-noncritical-neighbors},
respectively.
\end{lemma}
\begin{proof}
    We prove by contradiction. Suppose that at time $t$ there exist types $i$ and $j$ such that
$q_i(t)>0$, $q_j(t)>0$, and $(i,j)\in\mathcal M^g$ is a good match. Then we can pick a specific type-$i$ item and a specific type-$j$ item
that are \emph{available} at time $t$, i.e., they have arrived but have not been matched and have not abandoned by time $t$.
Let $t_i$ and $t_j$ denote their arrival times, respectively, and without loss of generality assume $t_i<t_j\le t$.

Consider the system at time $t_j$, when this type-$j$ item arrives. The chosen type-$i$ item is already present at time $t_j$
(because $t_i<t_j$) and remains available at time $t$; hence it is also available at time $t_j$.
Since $(i,j)$ is a good match and the type-$i$ item is available, the arriving type-$j$ item has at least one available good neighbor.
Algorithm~\ref{alg:shadow-matching-neighbor-set} therefore matches it on arrival, possibly with a different neighbor whose queue is longer. 
This contradicts the assumption that the type-$j$ item is available at time $t$.
Hence $q_i(t)q_j(t)=0$ for every $(i,j)\in\mathcal M^g$.

Now let $i\in\mathcal Q^+(t)$, so that $q_i(t)>0$. For every $j$ with $(i,j)\in\mathcal M^g$, the first part gives $q_j(t)=0$. Since $\mathcal N_i^+(t)$ and $\mathcal N_i^0(t)$,
defined in~\eqref{eq:fg-available-critical-neighbors}
and~\eqref{eq:fg-available-noncritical-neighbors}, respectively,
both require $q_j(t)>0$ and $(i,j)\in\mathcal M^g$,
they are empty.
\end{proof}

Building on this structural property, we now prove
Lemma~\ref{lmm:optimal_action}, which shows that $\bar{x}^b(t)$
is an optimal solution of problem~\eqref{eq:max_queue_prop}.

\begin{proof}[Proof of Lemma~\ref{lmm:optimal_action}]
Let $\delta(j) := \{i : M_{im}=M_{jm}=1\}$ denote the set of agent types that can be matched with type $j$. For any $x \in \mathcal{X}(t)$, we have
\begin{align}
\langle q(t), Mx\rangle_{\mathcal{I}^+} 
&= \sum_{i \in [n]} \sum_{j \in \delta(i) \cap \mathcal{I}^+} q_j(t) x_{(i,j)} \nonumber\\
&= \sum_{i \in [n] \setminus \mathcal{Q}^+(t)} \sum_{j \in \mathcal{N}_i^+(t)} q_j(t) x_{(i,j)}. \label{eq:objective_decomp}
\end{align}
The last equality holds due to  the following three observations:

\begin{enumerate}
    \item Since $q_j(t) = 0$ for $j \notin \mathcal{Q}^+(t)$, only terms with $j \in \mathcal{Q}^+(t)$ contribute.
    \item The constraint $x_m = 0$ for $m \notin \mathcal{M}^+(t)$ leaves only good matches, so the inner index is further restricted to those $j$ with $(i,j) \in \mathcal{M}^g$. Together with the previous observation, the surviving indices are exactly $\{ j \in \mathcal{I}^+ : q_j(t) > 0,\ (i,j) \in \mathcal{M}^g \} = \mathcal{N}_i^+(t)$ by~\eqref{eq:fg-available-critical-neighbors}.
    \item By Lemma~\ref{lmm:empty_neighbor}, $\mathcal{N}_i^+(t) = \emptyset$ for all $i \in \mathcal{Q}^+(t)$, eliminating all terms with $i \in \mathcal{Q}^+(t)$.
\end{enumerate}

From~\eqref{eq:objective_decomp}, since $i \in [n] \setminus \mathcal{Q}^+(t)$ and $j \in \mathcal{Q}^+(t)$, each variable $x_{(i,j)}$ appears exactly once. Thus, problem~\eqref{eq:max_queue_prop} decomposes into independent subproblems, one for each $i \in [n] \setminus \mathcal{Q}^+(t)$ with $\mathcal{N}_i^+(t) \neq \emptyset$, namely
\begin{equation}\label{eq:sub_max_queue}
\mathcal{P}_i(t):\quad
\begin{aligned}
    \max_{x_{(i,j)},\, j \in\mathcal{N}_i^+(t)} \quad &\sum_{j \in \mathcal{N}_i^+(t)} q_j(t) x_{(i,j)}\\
    \text{s.t.} \quad &\sum_{j \in \mathcal{N}_i^+(t)} x_{(i,j)} \leq \lambda^b_i, \quad x_{(i,j)} \geq 0 \quad \forall j \in \mathcal{N}_i^+(t).
\end{aligned}
\end{equation}
 We now show that $\bar{x}^b(t)=\sum_{i \in [n]} \lambda^b_i x^i(t)$ is feasible
for~\eqref{eq:max_queue_prop} and solves every subproblem. 

\textbf{Feasibility.} We first show that $\bar{x}^b(t) \in \mathcal{X}(t)$. By Algorithm~\ref{alg:shadow-matching-neighbor-set}, matches occur only for good match types and only when the arriving agent can be matched with a compatible agent currently in the queue. For any $m = (i,j) \notin \mathcal{M}^+(t)$, the algorithm never selects such a match, so $\bar{x}^b_m(t) = 0$.

Now fix $i \in [n] \setminus \mathcal{Q}^+(t)$, so that $q_i(t) = 0$; these are the only coordinates at which $\mathcal{X}(t)$ imposes a capacity constraint. For any $j \in \delta(i)$,
\begin{align*}
\bar{x}^b_{(i,j)}(t) = \lambda^b_i x^i_{(i,j)}(t) + \lambda^b_j x^j_{(i,j)}(t).
\end{align*}
When a type-$j$ agent arrives, no type-$i$ agent is waiting since $q_i(t) = 0$, so $x^j_{(i,j)}(t) = 0$. Therefore,
\begin{align}\label{eq:xbar-empty-type}
\bar{x}^b_{(i,j)}(t) = \lambda^b_i x^i_{(i,j)}(t) \quad \text{for all } j \in \delta(i),
\end{align}
and hence
\begin{align*}
    (M\bar{x}^b(t))_i = \sum_{j \in \delta(i)} \bar{x}^b_{(i,j)}(t) = \lambda^b_i \sum_{j \in \delta(i)} x^i_{(i,j)}(t) \leq \lambda^b_i,
\end{align*}
where the inequality holds since the algorithm selects at most one match per arrival. Thus $\bar{x}^b(t) \in \mathcal{X}(t)$.

\textbf{Optimality.} Consider $i \in [n] \setminus \mathcal{Q}^+(t)$ with
$\mathcal{N}^+_i(t) \neq \emptyset$. By~\eqref{eq:xbar-empty-type}, the
components of $\bar{x}^b(t)$ in problem~$\mathcal{P}_i(t)$ are
$\bar{x}^b_{(i,j)}(t) = \lambda^b_i x^i_{(i,j)}(t)$, $j \in \mathcal{N}_i^+(t)$.
Since the objective coefficients $q_j(t)$ of $\mathcal{P}_i(t)$ are
nonnegative and its only constraint is the total capacity $\lambda^b_i$,
$\mathcal{P}_i(t)$ is solved by giving the entire mass $\lambda^b_i$ to a
longest queue in $\mathcal{N}_i^+(t)$. This is exactly our decision rule:
Algorithm~\ref{alg:shadow-matching-neighbor-set} matches an arriving type-$i$
agent with the longest queue
$\arg\max_{j \in \mathcal{N}^+_i(t)} q_j(t)$,  so
$\bar{x}^b(t)$ gives mass $\lambda^b_i$ to this queue and zero to every other
type in $\mathcal{N}_i^+(t)$. Thus $\bar{x}^b(t)$ solves every
problem~$\mathcal{P}_i(t)$, and since $\bar{x}^b(t) \in \mathcal{X}(t)$, it
is optimal for~\eqref{eq:max_queue_prop}.

\end{proof}

\subsubsection{Proof of Lemma~\ref{lmm:construct_aux}}
\label{app:lmm:construct_aux}
Having proved Lemma~\ref{lmm:optimal_action}, we now turn to
Lemma~\ref{lmm:construct_aux}. Its proof uses the GPG condition
to construct a feasible auxiliary vector $\xi\in\mathcal{X}(t)$
satisfying~\eqref{eq:drift_construct}.

\begin{proof}[Proof of Lemma~\ref{lmm:construct_aux}]
Define
\[
\hat{\lambda}_i(t)=
\begin{cases}
    \lambda_i^b,
    & \text{if } i\notin\mathcal{I}^+\cap\mathcal{Q}^+(t),\\[4pt]
    \lambda_i^b+\epsilon/\sqrt n,
    & \text{if } i\in\mathcal{I}^+\cap\mathcal{Q}^+(t).
\end{cases}
\]
Since $|\mathcal{I}^+\cap\mathcal{Q}^+(t)|\leq n$, we have
\[
    \|\hat{\lambda}(t)-\lambda^b\|_2
    =\frac{\epsilon}{\sqrt n}
    \sqrt{|\mathcal{I}^+\cap\mathcal{Q}^+(t)|}
    \leq\epsilon.
\]
Thus, by Definition~\ref{def:GPG}, $U^*$ remains optimal
when $\lambda^b$ is replaced by $\hat{\lambda}(t)$.
Let $\hat{x}$ denote an optimal primal solution to \eqref{eq:reduced-fluid-primal} under arrival rates $\hat{\lambda}(t)$. We construct $\xi$ as follows:
\begin{align*}
\xi_m =
\begin{cases}
    \hat{x}_m, & \text{if } m \in \mathcal{M}^+(t),\\[4pt]
    0, & \text{otherwise}.
\end{cases}
\end{align*}

We first verify that $\xi \in \mathcal{X}(t)$. For $i \in [n] \setminus \mathcal{Q}^+(t)$,
\begin{align*}
[M\xi]_i  
\leq [M\hat{x}]_i  
\leq \hat{\lambda}_i(t)
= \lambda^b_i,
\end{align*}
confirming $\xi \in \mathcal{X}(t)$.

Next, we compute $[M\xi]_i$ for $i \in \mathcal{I}^+ \cap \mathcal{Q}^+(t)$:
\begin{align}
[M\xi]_i
= \sum_{j \in \delta(i)} \xi_{(i,j)}
= \sum_{j \in \delta(i)} \hat{x}_{(i,j)}
= \hat{\lambda}_i(t)
= \lambda_i^b+\frac{\epsilon}{\sqrt n}.
\label{eq:xi-critical-rate}
\end{align}
Since $U^*$ remains optimal when $\lambda^b$ is replaced by $\hat{\lambda}(t)$, and the dual constraints of~\eqref{eq:fluid-dual} do not involve the arrival rates, the set $\mathcal{M}^g$ is the same for both problems. The second equality then holds for $i \in \mathcal{Q}^+(t)$ as follows. If $(i,j) \in \mathcal{M}^g$, then since $i \in \mathcal{Q}^+(t)$ we have $(i,j) \in \mathcal{M}^+(t)$, so $\xi_{(i,j)} = \hat{x}_{(i,j)}$. If $(i,j) \notin \mathcal{M}^g$, then complementary slackness gives $\hat{x}_{(i,j)} = 0 = \xi_{(i,j)}$. The third equality holds because $i \in \mathcal{I}^+$ implies $U^*_i > 0$, which by complementary slackness forces $[M\hat{x}]_i = \hat{\lambda}_i(t)$. The last equality holds by the definition of $\hat{\lambda}(t)$, since $i \in \mathcal{I}^+ \cap \mathcal{Q}^+(t)$.
Finally, we establish the drift bound:
\begin{align*}
\langle q(t), \lambda^b - M\xi \rangle_{\mathcal{I}^+} 
&= -\frac{\epsilon}{\sqrt n}
    \sum_{i\in\mathcal{I}^+\cap\mathcal{Q}^+(t)}q_i(t)\\[3pt]
    &= -\frac{\epsilon}{\sqrt n}
    \|q(t)\|_{1,\mathcal{I}^+}\\[3pt]
    &\leq -\frac{\epsilon}{\sqrt n}
    \|q(t)\|_{\mathcal{I}^+}.
\end{align*}
The first equality holds since the restricted inner product ignores the coordinates outside $\mathcal{I}^+$ and $q_i(t) = 0$ for $i \notin \mathcal{Q}^+(t)$, together with~\eqref{eq:xi-critical-rate}. The second equality is the definition of $\|\cdot\|_{1,\mathcal{I}^+}$, again because $q_i(t) = 0$ off $\mathcal{Q}^+(t)$, and the final inequality holds since $\|z\|_{1,\mathcal{I}^+} \geq \|z\|_{\mathcal{I}^+}$. 
\end{proof}

\subsection{Proofs for Section~\ref{sec:step3}}
\label{app:step3}

We provide the proofs of the auxiliary lemmas used to
establish Proposition~\ref{prop:queue_length_bound}.

\subsubsection{Proof of Lemma~\ref{lmm:bounded_jump}}
\label{app:lmm:bounded_jump}
We begin with the bounded-jump property.
\begin{proof}[Proof of Lemma~\ref{lmm:bounded_jump}]
Since $\|\cdot\|_{\mathcal{I}^+}$ satisfies the triangle inequality, the reverse triangle inequality gives $\left| \|q(t+)\|_{\mathcal{I}^+} - \|q(t)\|_{\mathcal{I}^+} \right| \leq \|q(t+) - q(t)\|_{\mathcal{I}^+}$. It therefore suffices to bound $\|q(t+) - q(t)\|_{\mathcal{I}^+}$ by $1$ at each of the two event types identified in Lemma~\ref{lmm:generator}, an arrival and an abandonment. Throughout we use $\|e_k\|_{\mathcal{I}^+} = \mathbf{1}\{k \in \mathcal{I}^+\} \leq 1$.

Abandonment: If an agent of type $i$ abandons, the queue updates as $q(t+) = q(t) - e_i$. Consequently, $\|q(t+) - q(t)\|_{\mathcal{I}^+} = \|e_i\|_{\mathcal{I}^+} \leq 1$.

Arrival: Consider the arrival of a type-$i$ agent. Algorithm~\ref{alg:shadow-matching-neighbor-set} either attempts a match or leaves the agent in the queue.
\begin{enumerate}
    \item No Match: If no match is executed, $q(t+) = q(t) + e_i$. Hence, $\|q(t+) - q(t)\|_{\mathcal{I}^+} = \|e_i\|_{\mathcal{I}^+} \leq 1$.
    \item Match Executed: The selected match involves the arriving type $i$ and a waiting type $j \neq i$, so $Mx^i(t) = e_i + e_j$ and $q(t+) = q(t) + e_i - (e_i + e_j) = q(t) - e_j$. Hence $\|q(t+) - q(t)\|_{\mathcal{I}^+} = \|e_j\|_{\mathcal{I}^+} \leq 1$.
\end{enumerate}
In all cases, the jump size is bounded by 1.
\end{proof}

\subsubsection{Proof of Lemma~\ref{lmm:queue_negdrift}}
\label{app:lmm:queue_negdrift}
Next, we use the negative drift bound for
$\|q(t)\|_{\mathcal I^+}^2$ in
Proposition~\ref{prop:queue_square_negdrift} to derive
the drift bounds for $\|q(t)\|_{\mathcal I^+}$.

\begin{proof}[Proof of Lemma~\ref{lmm:queue_negdrift}]
Applying the generator~\eqref{eq:generator} with $f(q) = \|q\|_{\mathcal{I}^+}$ yields
\begin{align}
    \mathcal{L}\|q(t)\|_{\mathcal{I}^+} &= \sum_{i=1}^n \lambda_i^s \bigl[ \|q(t) + e_i - Mx^i(t)\|_{\mathcal{I}^+} - \|q(t)\|_{\mathcal{I}^+} \bigr] + \sum_{i \in \mathcal{I}^+} \theta_i q_i(t) \bigl[ \|q(t) - e_i\|_{\mathcal{I}^+} - \|q(t)\|_{\mathcal{I}^+} \bigr], \label{eq:generator_norm_expansion}
\end{align}
where the abandonment sum is restricted to $i \in \mathcal{I}^+$ since abandonments from $i \notin \mathcal{I}^+$ do not affect $\|q(t)\|_{\mathcal{I}^+}$. By Lemma~\ref{lmm:bounded_jump}, each jump satisfies $| \|q(t^+)\|_{\mathcal{I}^+} - \|q(t)\|_{\mathcal{I}^+} | \leq 1$, so
\begin{align}
    \mathcal{L}\|q(t)\|_{\mathcal{I}^+} \leq s + \|\theta\|_{\mathcal{I}^+} \|q(t)\|_{\mathcal{I}^+}, \label{eq:generator_upper}
\end{align}
where we use $\sum_{i=1}^n \lambda_i^s  = s\|\lambda^b\|_1 = s$ and the Cauchy–Schwarz inequality.

To obtain negative drift for large queue lengths, we use the inequality $\|a\| - \|b\| \leq \frac{\|a\|^2 - \|b\|^2}{2\|b\|}$, valid for any $b \neq 0$. Applying this to each term in~\eqref{eq:generator_norm_expansion}, we have
\begin{align*}
    \mathcal{L}\|q(t)\|_{\mathcal{I}^+} &\leq \sum_{i=1}^n \lambda_i^s \frac{\|q(t) + e_i - Mx^i(t)\|_{\mathcal{I}^+}^2 - \|q(t)\|_{\mathcal{I}^+}^2}{2\|q(t)\|_{\mathcal{I}^+}} + \sum_{i \in \mathcal{I}^+} \theta_i q_i(t) \frac{\|q(t) - e_i\|_{\mathcal{I}^+}^2 - \|q(t)\|_{\mathcal{I}^+}^2}{2\|q(t)\|_{\mathcal{I}^+}} \notag \\
    &= \frac{\mathcal{L}\|q(t)\|_{\mathcal{I}^+}^2}{2\|q(t)\|_{\mathcal{I}^+}}, 
\end{align*}
where the equality holds since the numerator is the generator~\eqref{eq:generator} with $f(q) = \|q\|_{\mathcal{I}^+}^2$.
By Proposition~\ref{prop:queue_square_negdrift}, $\mathcal{L}\|q(t)\|_{\mathcal{I}^+}^2 \leq 2s ( -\frac{\varepsilon}{\sqrt{n}} \|q(t)\|_{\mathcal{I}^+} + \frac{1}{2} )$. Substituting yields
\begin{align}
    \mathcal{L}\|q(t)\|_{\mathcal{I}^+} \leq s \left( -\frac{\varepsilon}{\sqrt{n}} + \frac{1}{2\|q(t)\|_{\mathcal{I}^+}} \right). \label{eq:generator_neg}
\end{align}
When $\|q(t)\|_{\mathcal{I}^+} \geq \frac{\sqrt{n}}{\varepsilon}$, the second term in~\eqref{eq:generator_neg} is at most $\frac{\varepsilon}{2\sqrt{n}}$, yielding $\mathcal{L}\|q(t)\|_{\mathcal{I}^+} \leq -\frac{s\varepsilon}{2\sqrt{n}}$. When $\|q(t)\|_{\mathcal{I}^+} < \frac{\sqrt{n}}{\varepsilon}$, substituting into~\eqref{eq:generator_upper} gives $\mathcal{L}\|q(t)\|_{\mathcal{I}^+} \leq s + \frac{\sqrt{n}\|\theta\|_{\mathcal{I}^+}}{\varepsilon}$.

\end{proof}

\subsubsection{Proof of Lemma~\ref{lmm:general_drift}}
\label{app:lmm:general_drift}
Finally, we prove Lemma~\ref{lmm:general_drift}, which
establishes a uniform bound on the expected value of
a function under the stated drift and jump conditions.

\begin{proof}[Proof of Lemma~\ref{lmm:general_drift}]
Consider the exponential Lyapunov function $V(y) = e^{\alpha f(y)}$, where $\alpha > 0$ is a parameter to be determined later. By condition~\ref{cond:jump-rep} applied to $V$,
\begin{align}
    \mathcal{L}V(y) = \sum_{j \in \mathcal{J}} \lambda_j(y) \left( e^{\alpha (f(y) + \Delta_j f(y))} - e^{\alpha f(y)} \right) = V(y) \sum_{j \in \mathcal{J}} \lambda_j(y) \left( e^{\alpha \Delta_j f(y)} - 1 \right). \label{eq:exp_generator}
\end{align}
We partition the jumps into negative jumps $\mathcal{J}^-(y) = \{j : \Delta_j f(y) < 0\}$ and positive jumps $\mathcal{J}^+(y)$. We then bound the exponential term $e^{\alpha \Delta_j f(y)} - 1$ for each case separately.

For $j \in \mathcal{J}^-(y)$, let $z = \Delta_j f(y) \in [-K, 0)$. Using the inequality $e^u \leq 1 + u + \frac{u^2}{2}$ for $u \leq 0$, and observing that $z^2 \leq -Kz$ on this interval, we have
\begin{align}
    e^{\alpha z} - 1 \leq \alpha z + \frac{\alpha^2 z^2}{2} \leq \alpha z - \frac{\alpha^2 K z}{2} = \alpha \left(1 - \frac{\alpha K}{2}\right) z. \label{eq:neg-jump-bound}
\end{align}
For $j \in \mathcal{J}^+(y)$, let $z = \Delta_j f(y) \in (0, K]$. Restricting $\alpha \in (0, 1/K]$ ensures $\alpha z \in (0, 1]$. Using the inequality $e^u - 1 \leq u + (e-2)u^2$ for $u \in [0,1]$, we obtain
\begin{align}
    e^{\alpha z} - 1 \leq \alpha z + (e - 2) \alpha^2 z^2. \label{eq:pos-jump-bound}
\end{align}
Substituting~\eqref{eq:neg-jump-bound} and~\eqref{eq:pos-jump-bound} into~\eqref{eq:exp_generator} yields
\begin{align}
    \frac{\mathcal{L}V(y)}{V(y)} &\leq \sum_{j \in \mathcal{J}^-(y)} \lambda_j(y) \cdot \alpha \left(1 - \frac{\alpha K}{2}\right) \Delta_j f(y) + \sum_{j \in \mathcal{J}^+(y)} \lambda_j(y) \left( \alpha \Delta_j f(y) + (e - 2) \alpha^2 (\Delta_j f(y))^2 \right) \notag \\
    &= \alpha \left(1 - \frac{\alpha K}{2}\right) \mathcal{L}f(y) + \frac{\alpha^2 K}{2} \sum_{j \in \mathcal{J}^+(y)} \lambda_j(y) \Delta_j f(y) + (e - 2) \alpha^2 \sum_{j \in \mathcal{J}^+(y)} \lambda_j(y) (\Delta_j f(y))^2 \notag \\
    &\leq \alpha \left(1 - \frac{\alpha K}{2}\right) \mathcal{L}f(y) + (e - \frac{3}{2}) \Lambda K^2  \alpha^2, \label{eq:exp-generator-bound}
\end{align}
where the last inequality uses conditions~\ref{cond:bounded-jump} and~\ref{cond:bounded-rate}, specifically $\sum_{j \in \mathcal{J}^+(y)} \lambda_j(y) \Delta_j f(y) \leq K \Lambda$ and $\sum_{j \in \mathcal{J}^+(y)} \lambda_j(y) (\Delta_j f(y))^2 \leq K^2 \Lambda$.

We now select $\alpha$ to establish a global drift condition. In the region $\{y : f(y) \geq D\}$, condition~\ref{cond:drift} gives $\mathcal{L}f(y) \leq -\delta$, so~\eqref{eq:exp-generator-bound} becomes
\begin{align*}
    \frac{\mathcal{L}V(y)}{V(y)} \leq -\alpha \delta + \alpha^2 C_1 ,
\end{align*}
where $C_1 = \frac{K\delta}{2} + (e - \frac{3}{2})\Lambda K^2 $. 
Choosing $\alpha^* = \min( \frac{1}{K}, \frac{\delta}{2C_1} )$ ensures 
\begin{align}
    \mathcal{L}V(y) \leq -\gamma V(y)\label{eq:exp_large_drift}
\end{align}
with $\gamma = \alpha^* \delta - (\alpha^*)^2 C_1 > 0$.

In the region $\{y : f(y) < D\}$, condition~\ref{cond:drift} gives $\mathcal{L}f(y) \leq B$, so~\eqref{eq:exp-generator-bound} becomes
\begin{align}
    \frac{\mathcal{L}V(y)}{V(y)} &\leq \alpha^* B + \left( -\frac{KB}{2} + (e - 1.5)\Lambda K^2 \right) (\alpha^*)^2 \notag\\
    &\leq \alpha^* B + C_1 (\alpha^*)^2 \notag\\
    &= - \gamma + \alpha^*(B + \delta) ,\label{eq:exp_small_drift}
\end{align}
where the second inequality follows because the coefficient of $(\alpha^*)^2$ is less than $C_1$, and the final equality uses the definition $\gamma = \alpha^* \delta - (\alpha^*)^2 C_1$.
Since $V(y) \leq e^{\alpha^* D}$ in this region, combining~\eqref{eq:exp_large_drift} and~\eqref{eq:exp_small_drift} yields the global drift condition
\begin{align*}
    \mathcal{L}V(y) \leq -\gamma V(y) + C_{\mathrm{const}}, \quad \text{where } C_{\mathrm{const}} = \alpha^*(B + \delta) e^{\alpha^* D}.
\end{align*}

Applying Dynkin's formula, we obtain
\begin{align*}
    \frac{d}{dt} \mathbb{E}[V(Y_t)] = \mathbb{E}[\mathcal{L}V(Y_t)].
\end{align*}
Substituting the drift bound $\mathcal{L}V(y) \leq -\gamma V(y) + C_{\mathrm{const}}$ into the right-hand side yields the differential inequality
\begin{align*}
    \frac{d}{dt} \mathbb{E}[V(Y_t)] \leq -\gamma \mathbb{E}[V(Y_t)] + C_{\mathrm{const}}.
\end{align*}
By Grönwall's inequality, we obtain
\begin{align*}
    \mathbb{E}[V(Y_t)] &\leq e^{-\gamma t} \mathbb{E}[V(Y_0)] + \frac{C_{\mathrm{const}}}{\gamma} (1 - e^{-\gamma t}) \\
    &\leq \max\left\{ \mathbb{E}[V(Y_0)], \frac{C_{\mathrm{const}}}{\gamma} \right\},
\end{align*}
where the second inequality follows because the expression is a convex combination with weights summing to 1.
Finally, applying Jensen's inequality implies 
\begin{align*}
    e^{\alpha^* \mathbb{E}[f(Y_t)]} \leq \mathbb{E}[e^{\alpha^* f(Y_t)}] = \mathbb{E}[V(Y_t)].
\end{align*}
 Taking logarithms and substituting the constant $\frac{C_{\mathrm{const}}}{\gamma} = e^{\alpha^* D} \left( \frac{B + \delta}{\delta - \alpha^* C_1} \right)$ yields
\begin{align*}
    \mathbb{E}[f(Y_t)] \leq \frac{1}{\alpha^*} \max\left\{ \alpha^* f(Y_0),\, \ln \left( \frac{B + \delta}{\delta - \alpha^* C_1} \right) + \alpha^* D \right\}.
\end{align*}
\end{proof}

\subsection{Bounded Dual Solutions}
\label{app:dual_bound}

We state and prove Lemma~\ref{lmm:dual_bound}, which bounds an
optimal dual solution of the fluid LP.

\begin{lemma}
\label{lmm:dual_bound}
Let \(U^*\) be any optimal solution of~\eqref{eq:fluid-dual}
at arrival rates \(\lambda^b\). Then
\begin{equation*}
    \|U^*\|_\infty\leq\widetilde r_{\max},
    \qquad
    \|U^*\|\leq\sqrt n\,\widetilde r_{\max}.
\end{equation*}
\end{lemma}
\begin{proof}[Proof of Lemma~\ref{lmm:dual_bound}]
Fix \(k\in\mathcal N\), and suppose that
\(U_k^*>\widetilde r_{\max}\). Define
\[
    U'=U^*-(U_k^*-\widetilde r_{\max})e_k.
\]
For every match \(m=(k,j)\), we have
\[
    U'_k+U'_j=\widetilde r_{\max}+U_j^*
    \geq\widetilde r_m.
\]
All other dual constraints are unchanged, and \(U'\geq0\), so \(U'\)
is dual feasible. However,
\[
    (\lambda^b)^\top U'-(\lambda^b)^\top U^*
    =-\lambda_k^b(U_k^*-\widetilde r_{\max})<0,
\]
contradicting the optimality of \(U^*\). Since \(k\) is arbitrary,
\[
    \|U^*\|_\infty\leq\widetilde r_{\max},
    \qquad
    \|U^*\|^2\leq n\|U^*\|_\infty^2
    \leq n\widetilde r_{\max}^2.
\]
\end{proof}

\section{Proof of Proposition~\ref{prop:no_gpg}}
\label{app:proof_no_gpg}

In this section, we prove the regret bound in
Proposition~\ref{prop:no_gpg} without the GPG condition. We first establish
an \(O(\sqrt{s})\) bound on \(\mathbb E[\|q(t)\|_{\mathcal I^+}]\),
and then use this bound to derive the regret guarantee.

Under the GPG condition,  the proof of Proposition~\ref{prop:queue_square_negdrift} shows that
\[
    \langle q(t),\lambda^b-M\bar{x}^b(t)\rangle_{\mathcal I^+}
    \leq-\frac{\varepsilon}{\sqrt{n}}\|q(t)\|_{\mathcal I^+},
\]
which provides a downward drift contribution for
\(\|q(t)\|_{\mathcal I^+}^2\). Without the GPG condition, we obtain the
weaker bound
\[
    \langle q(t),\lambda^b-M\bar{x}^b(t)\rangle_{\mathcal I^+}\leq0.
\]
Thus, this term remains nonpositive and does not contribute to the growth
of \(\|q(t)\|_{\mathcal I^+}^2\). Following a proof strategy similar to
that of Lemma~\ref{lmm:construct_aux}, we construct an auxiliary
vector \(\xi(t)\in\mathcal X(t)\) to establish this bound in
Lemma~\ref{lmm:no-gpg-matching}.

\begin{lemma}
\label{lmm:no-gpg-matching}
Under the longest-queue policy described in
Algorithm~\ref{alg:shadow-matching-neighbor-set}, given the queue state
\(q(t)\), we have
\begin{equation*}
    \langle q(t),\lambda^b-M\bar{x}^b(t)\rangle_{\mathcal I^+}
    \leq0.
\end{equation*}
\end{lemma}

\begin{proof}
Let \(x^*\) be an optimal solution of~\eqref{eq:reduced-fluid-primal}
at rates \(\lambda^b\). Complementary slackness gives
\[
    x_m^*=0\quad\text{for }m\notin\mathcal M^g,
    \qquad
    (Mx^*)_i=\lambda_i^b\quad\text{for }i\in\mathcal I^+.
\]
Define
\[
    \xi_m(t)=
    \begin{cases}
        x_m^*,&m\in\mathcal M^+(t),\\
        0,&m\notin\mathcal M^+(t).
    \end{cases}
\]
For every \(i\notin\mathcal Q^+(t)\),
\[
    (M\xi(t))_i\leq(Mx^*)_i\leq\lambda_i^b.
\]
Together with the definition of \(\xi(t)\), this shows that
\(\xi(t)\in\mathcal X(t)\), where \(\mathcal X(t)\) is defined in~\eqref{eq:auxiliary-feasible-set}.

Let \(\delta(i)\) denote the set of agent types that can be matched with
type \(i\). For \(i\in\mathcal I^+\cap\mathcal Q^+(t)\), we have
\[
    (Mx^*)_i
    =\sum_{\substack{j\in\delta(i)\\(i,j)\in\mathcal M^g}}x^*_{(i,j)}
    =\sum_{\substack{j\in\delta(i)\\(i,j)\in\mathcal M^g}}\xi_{(i,j)}(t)
    =(M\xi(t))_i,
\]
where the first equality follows from \(x^*_{(i,j)}=0\) for
\((i,j)\notin\mathcal M^g\); the second holds because
\((i,j)\in\mathcal M^g\) and \(i\in\mathcal Q^+(t)\) imply
\((i,j)\in\mathcal M^+(t)\), so \(\xi_{(i,j)}(t)=x^*_{(i,j)}\)
by definition; and the last follows from \(\xi_{(i,j)}(t)=0\) for
\((i,j)\notin\mathcal M^g\).

Thus, \((M\xi(t))_i=(Mx^*)_i=\lambda_i^b\) for
\(i\in\mathcal I^+\cap\mathcal Q^+(t)\). For
\(i\in\mathcal I^+\setminus\mathcal Q^+(t)\), we have \(q_i(t)=0\).
Therefore,
\[
    \langle q(t),\lambda^b-M\xi(t)\rangle_{\mathcal I^+}=0.
\]
By Lemma~\ref{lmm:optimal_action}, \(\bar{x}^b(t)\) maximizes
\(\langle q(t),Mx\rangle_{\mathcal I^+}\) over \(\mathcal X(t)\).
Since \(\xi(t)\in\mathcal X(t)\),
\[
    \langle q(t),\lambda^b-M\bar{x}^b(t)\rangle_{\mathcal I^+}
    \leq
    \langle q(t),\lambda^b-M\xi(t)\rangle_{\mathcal I^+}
    =0.
\]
\end{proof}

Next, we combine Lemma~\ref{lmm:no-gpg-matching} with the downward drift
induced by abandonment to establish an \(O(\sqrt{s})\) bound on
\(\mathbb E[\|q(t)\|_{\mathcal I^+}]\) in
Lemma~\ref{lmm:no-gpg-queues}.

\begin{lemma}
\label{lmm:no-gpg-queues}
Let \(\theta_{\min}=\min_{i\in\mathcal N}\theta_i>0\).
For every \(t\geq0\),
\begin{equation*}
    \mathbb E\bigl[\|q(t)\|_{\mathcal I^+}^2\bigr]
    \leq\frac{s}{\theta_{\min}}
        \bigl(1-e^{-\theta_{\min}t}\bigr).
\end{equation*}
Consequently,
\begin{equation*}
    \sup_{t\geq0}\mathbb E\bigl[\|q(t)\|_{\mathcal I^+}\bigr]
    \leq\sqrt{\frac{s}{\theta_{\min}}}.
\end{equation*}
\end{lemma}

\begin{proof}
Define \(V(q(t))=\|q(t)\|_{\mathcal I^+}^2\).
Combining the expansion~\eqref{eq:drift-expanded-abandonment} in the
proof of Lemma~\ref{lmm:generator_bound} with the
bound~\eqref{eq:jump-bound}, and using \(\lambda^s=s\lambda^b\) and
\(\bar{x}(t)=s\bar{x}^b(t)\), gives
\begin{align*}
    \mathcal L V(q(t))
    &\leq
    2s\langle q(t),\lambda^b-M\bar{x}^b(t)\rangle_{\mathcal I^+}
    +s 
    +
    \sum_{i\in\mathcal I^+}\theta_iq_i(t)\bigl[-2q_i(t)+1\bigr].
\end{align*}
The first term is nonpositive by Lemma~\ref{lmm:no-gpg-matching}.
Since \(q_i(t)\) is a nonnegative integer,
\(q_i(t)\bigl[-2q_i(t)+1\bigr]\leq-q_i(t)^2\). Hence
\begin{equation}
    \mathcal L V(q(t))
    \leq s-\sum_{i\in\mathcal I^+}\theta_iq_i(t)^2
    \leq s-\theta_{\min}V(q(t)).
    \label{eq:no-gpg-drift}
\end{equation}

Applying Dynkin's formula gives
\[
    \frac{\diff}{\diff t}\mathbb E[V(q(t))]
    =\mathbb E[\mathcal L V(q(t))].
\]
Substituting the drift bound~\eqref{eq:no-gpg-drift} yields
\[
    \frac{\diff}{\diff t}\mathbb E[V(q(t))]
    \leq s-\theta_{\min}\mathbb E[V(q(t))].
\]
Since \(q(0)=0\), we have \(\mathbb E[V(q(0))]=0\), and
Gr\"onwall's inequality gives
\[
    \mathbb E[V(q(t))]
    \leq\frac{s}{\theta_{\min}}
        \bigl(1-e^{-\theta_{\min}t}\bigr).
\]
By Cauchy--Schwarz, we obtain
\[
    \mathbb E\bigl[\|q(t)\|_{\mathcal I^+}\bigr]
    \leq\sqrt{\mathbb E\bigl[\|q(t)\|_{\mathcal I^+}^2\bigr]}
    \leq\sqrt{\frac{s}{\theta_{\min}}},
    \qquad t\geq0.
\]
\end{proof}

Having already expressed the regret in terms of the queue lengths in
Proposition~\ref{prop:regret_decomposition}, we now combine this representation
with the queue-length bound in Lemma~\ref{lmm:no-gpg-queues} to prove
Proposition~\ref{prop:no_gpg}.

\begin{proof}[Proof of Proposition~\ref{prop:no_gpg}]
Applying the regret representation in Proposition~\ref{prop:regret_decomposition}
at arrival rates \(s\lambda^b\), we obtain
\begin{align*}
    \mathcal R_{T,s}^{\pi}
    &={}
    \frac{1}{T}\langle U^*,\mathbb E[q(T)]\rangle_{\mathcal I^+}
    +\frac{1}{T}\int_0^T
        \mathbb E\bigl[\langle U^*,\theta q(t)\rangle_{\mathcal I^+}\bigr]
        \,\diff t\\
    &\leq
    \frac{\sqrt n\,\widetilde r_{\max}}{T}
        \mathbb E\bigl[\|q(T)\|_{\mathcal I^+}\bigr]
    +\frac{\|\theta\|_\infty\sqrt n\,\widetilde r_{\max}}{T}
        \int_0^T\mathbb E\bigl[\|q(t)\|_{\mathcal I^+}\bigr]
        \,\diff t\\
    &\leq
    \left(\frac{1}{T}+\|\theta\|_\infty\right)
    \sqrt n\,\widetilde r_{\max}\sqrt{\frac{s}{\theta_{\min}}}.
\end{align*}
The first inequality follows from Cauchy--Schwarz and the bound
\(\|U^*\|\leq\sqrt n\,\widetilde r_{\max}\) in
Lemma~\ref{lmm:dual_bound}. The last inequality follows from the
queue-length bound in Lemma~\ref{lmm:no-gpg-queues},
proving~\eqref{eq:no-gpg-regret}.
Finally, \(F_s^*=sF^*(\lambda^b)\), so
\[
    \rho_{T,s}^{\pi}
    =1-\frac{\mathcal R_{T,s}^{\pi}}{sF^*(\lambda^b)}
    \geq
    1-
    \frac{(T^{-1}+\|\theta\|_\infty)\sqrt n\,\widetilde r_{\max}}
         {F^*(\lambda^b)\sqrt{\theta_{\min}}}\,
    \frac{1}{\sqrt{s}}.
\]
For \(T\geq1\), this is at least \(1-C/\sqrt{s}\), where \(C\) is
independent of \(T\) and \(s\), proving~\eqref{eq:no-gpg-ratio}.
\end{proof}
\end{APPENDICES}
\end{document}